\documentclass[aop,preprint]{imsart}

\setattribute{journal}{name}{}

\usepackage[utf8]{inputenc}
\usepackage[T1]{fontenc}
\usepackage{totpages}
\usepackage{graphicx}
\usepackage{amsmath, amsthm, amssymb}
\usepackage{upgreek}
\usepackage{mathrsfs}
\usepackage{pdfsync}
\usepackage[a4paper,top=3cm,bottom=3cm,left=2.5cm,right=2.5cm, bindingoffset=0mm]{geometry}

\usepackage[usenames,dvipsnames]{xcolor}
\usepackage[inline,shortlabels]{enumitem}
\usepackage[hyphens]{url}									%per gli url in bliografia
\usepackage{verbatim}								%per commentare cose lunghe
\usepackage{dsfont}									%doublestroke \mathds

\usepackage{mathtools}
\usepackage{xparse}
\newcommand{\texorpdfstring}[2]{#1}

\usepackage{mathtools}

\allowdisplaybreaks

\arxiv{}

\startlocaldefs

\usepackage{xparse}

\ExplSyntaxOn
\NewDocumentCommand{\makeabbrev}{mmm}
 {
  \yoruk_makeabbrev:nnn { #1 } { #2 } { #3 }
 }

\cs_new_protected:Npn \yoruk_makeabbrev:nnn #1 #2 #3
 {
  \clist_map_inline:nn { #3 }
   {
    \cs_new_protected:cpn { #2 } { #1 { ##1 } }
   }
 }
 \ExplSyntaxOff

\makeabbrev{\textbf}{tbf#1}{a,b,c,d,e,f,g,h,i,j,k,l,m,n,o,p,q,r,s,t,u,v,w,x,y,z,A,B,C,D,E,F,G,H,I,J,K,L,M,N,O,P,Q,R,S,T,U,V,W,X,Y,Z}

\makeabbrev{\textbf}{bf#1}{a,b,c,d,e,f,g,h,i,j,k,l,m,n,o,p,q,r,s,t,u,v,w,x,y,z,A,B,C,D,E,F,G,H,I,J,K,L,M,N,O,P,Q,R,S,T,U,V,W,X,Y,Z}

\makeabbrev{\textsf}{tsf#1}{a,b,c,d,e,f,g,h,i,j,k,l,m,n,o,p,q,r,s,t,u,v,w,x,y,z,A,B,C,D,E,F,G,H,I,J,K,L,M,N,O,P,Q,R,S,T,U,V,W,X,Y,Z}

\makeabbrev{\mathsf}{mss#1}{a,b,c,d,e,f,g,h,i,j,k,l,m,n,o,p,q,r,s,t,u,v,w,x,y,z,A,B,C,D,E,F,G,H,I,J,K,L,M,N,O,P,Q,R,S,T,U,V,W,X,Y,Z}

\makeabbrev{\mathfrak}{mf#1}{a,b,c,d,e,f,g,h,i,j,k,l,m,n,o,p,q,r,s,t,u,v,w,x,y,z,A,B,C,D,E,F,G,H,I,J,K,L,M,N,O,P,Q,R,S,T,U,V,W,X,Y,Z}

\makeabbrev{\mathrm}{mrm#1}{a,b,c,d,e,f,g,h,i,j,k,l,m,n,o,p,q,r,s,t,u,v,w,x,y,z,A,B,C,D,E,F,G,H,I,J,K,L,M,N,O,P,Q,R,S,T,U,V,W,X,Y,Z}

\makeabbrev{\mathbf}{mbf#1}{a,b,c,d,e,f,g,h,i,j,k,l,m,n,o,p,q,r,s,t,u,v,w,x,y,z,A,B,C,D,E,F,G,H,I,J,K,L,M,N,O,P,Q,R,S,T,U,V,W,X,Y,Z}

\makeabbrev{\mathcal}{mc#1}{A,B,C,D,E,F,G,H,I,J,K,L,M,N,O,P,Q,R,S,T,U,V,W,X,Y,Z}

\makeabbrev{\mathbb}{mbb#1}{A,B,C,D,E,F,G,H,I,J,K,L,M,N,O,P,Q,R,S,T,U,V,W,X,Y,Z}

\makeabbrev{\mathscr}{ms#1}{A,B,C,D,E,F,G,H,I,J,K,L,M,N,O,P,Q,R,S,T,U,V,W,X,Y,Z}

\makeabbrev{\mathrm}{#1}{
Id,id,ran,rk,diag,stab,ann,conv,pr,ev,tr,End,Hom,sgn,im,op,can,fin,ext,red,tot,Leb,
rot,usc,lsc,Lip,lip,bSymLip,osc,AC,loc,coz,z,
supp,Opt,Adm,Cpl,Geo,GeoOpt,GeoAdm,GeoCpl,reg,res,graph,
bd,co,Ric,Exp,dExp,dist,seg,Seg,cut,fcut,Cut,SDiff,Iso,Isom,diam,cl,Homeo,Diff,Der,vol,dvol,inj,relint, Graph, sub,
var,law,Var,Poi,Gam,pa,so,iso,fs,inv,pqi,mix,erg,
TestF,
}

\makeabbrev{\mathsf}{#1}{CD,BE,MCP,Ent,wMTW,MTW,Ch,RCD,EVI,Rad,dRad,SL,cSL,dSL,ScL,Irr,SC,wFe,VA}

\makeabbrev{\mathsc}{#1}{mmaf,cg}

\newcommand{\T}{\mathfrak{T}} %TOPOLOGY
\newcommand{\eps}{\varepsilon}

\renewcommand{\complement}{\mathrm{c}}

\newcommand{\mathsc}[1]{\text{\textsc{#1}}}
\newcommand{\emparg}{{\,\cdot\,}}
\newcommand{\dint}[2][]{\;\sideset{^{\scriptstyle{#1}}\!\!\!\!}{_{#2}^{\scriptscriptstyle\oplus}}\int}
\newcommand{\forallae}[1]{{\textrm{\,for ${#1}$-a.e.~}}}

\renewcommand{\Cap}{\mathrm{cap}}
\newcommand{\dom}[1]{\msD(#1)}

\DeclareMathOperator{\eqdef}{\coloneqq}

\let\epsilon\varepsilon

\newcommand{\longrar}{\longrightarrow}
\newcommand{\rar}{\rightarrow}

\newcommand{\diff}{\mathop{}\!\mathrm{d}}						%Differenziale esatto

\newcommand{\abs}[1]{\left\lvert#1\right\rvert}						%Modulo
\newcommand{\norm}[1]{\left\lVert#1\right\rVert}					%Norma
\newcommand{\set}[1]{\left\{#1\right\}}							%Insieme, graffe
\newcommand{\ttset}[1]{\{#1\}}									%Insieme, graffe
\newcommand{\tonde}[1]{\left(#1\right)}							%Tonde
\newcommand{\ttonde}[1]{\big({#1}\big)}
\newcommand{\class}[2][]{\left[#2\right]_{#1}}						%Measure classes
\newcommand{\sym}[1]{{\scriptscriptstyle{(#1)}}}

\newcommand{\rep}[1]{\hat{#1}}								%Rappresentante
\newcommand{\scalar}[2]{\left\langle #1 \,\middle |\, #2\right\rangle}		%Prodotto scalare	
\newcommand{\seq}[1]{\tonde{#1}}								%Successione
\newcommand{\tseq}[1]{{\big(#1\big)}}
\newcommand{\ttseq}[1]{{(#1)}}
\newcommand{\Cz}{\mcC_0}									%funzioni continue e evanescenti

\newcommand{\pfwd}{\sharp}
\DeclareMathOperator*{\esssup}{esssup}

\DeclareMathOperator{\car}{\mathbf 1}

\newcommand{\R}{{\mathbb R}}

\newcommand{\restr}{\big\lvert}
\newcommand{\mrestr}[1]{\!\!\downharpoonright_{#1}}

\newcommand{\iref}[1]{\ref{#1}}

\newcommand{\comma}{\,\mathrm{,}\;\,}

\newcommand{\fstop}{\,\mathrm{.}}

\DeclareMathOperator{\zero}{{\mathbf 0}}

\usepackage{scrextend}						%KOMA: confligge quindi deve stare qui

\let\temp\phi
\let\phi\varphi
\let\varphi\temp

\numberwithin{equation}{section}
\theoremstyle{plain}
\newtheorem{theorem}{Theorem}[section]
\newtheorem{proposition}[theorem]{Proposition}%[section]
\newtheorem{lemma}[theorem]{Lemma}%[section]
\newtheorem{corollary}[theorem]{Corollary}%[section]

\theoremstyle{definition}
\newtheorem{definition}[theorem]{Definition}%[section]
\newtheorem*{defs*}{Definition}%[section]
\theoremstyle{remark}
\newtheorem{remark}[theorem]{Remark}%[section]
\newtheorem{example}[theorem]{Example}%[section]
\newtheorem{Assumption}[theorem]{Assumption}%[section]

\newcommand{\purple}[1]{{{#1}}}

\endlocaldefs

\begin{document}

\begin{frontmatter}
\title{Ergodic Decompositions of Dirichlet forms\\
under order isomorphisms%\thanksref{t0}
}
\runtitle{Ergodic Decompositions of Dirichlet forms under order isomorphisms}
%\thankstext{T1}{Footnote to the title with the ``thankstext'' command.}

\begin{aug}
\author{\fnms{Lorenzo} \snm{Dello Schiavo}\thanksref{t1}\ead[label=e1]{lorenzo.delloschiavo@ist.ac.at}}
and
\author{\fnms{Melchior} \snm{Wirth}\thanksref{t1}\ead[label=e2]{melchior.wirth@ist.ac.at}}

\runauthor{L.~Dello~Schiavo and M.~Wirth}

\affiliation{Institute of Science and Technology Austria}

\address{IST Austria\\
Am Campus 1\\
3400 Klosterneuburg\\
Austria\\
\printead{e1}\\
\printead{e2}
%\phantom{E-mail:\ }\printead*{e2}
}

%\thankstext{t0}{\purple{Acknowledgements to people}}
\thankstext{t1}{Research supported by the Austrian Science Fund (FWF) grant F65 at the Institute of Science and Technology Austria and by the European Research Council (ERC) (grant agreement No.~716117 awarded to Prof.~Dr.~Jan Maas).
%\\
L.D.S.~gratefully acknowledges funding of his current position by the Austrian Science Fund (FWF) through the ESPRIT Programme (grant No.~208).
%\\
M.W.~gratefully acknowledges funding of his current position by the Austrian Science Fund (FWF) through the ESPRIT Programme (grant No.~156).
}

\end{aug}

\begin{abstract}
We study ergodic decompositions of Dirichlet spaces under intertwining via unitary order isomorphisms.
We show that the ergodic decomposition of a quasi-regular Dirichlet space is unique up to a unique isomorphism of the indexing space.
Furthermore, every unitary order isomorphism intertwining two quasi-regular Dirichlet spaces is decomposable over their ergodic decompositions up to conjugation via an isomorphism of the corresponding indexing spaces.
\end{abstract}

\vspace{.5cm}
\today
\vspace{.5cm}



\end{frontmatter}

%\tableofcontents

\section{Introduction}

How much geometric information can be recovered from associated analytic or probabilistic objects is the underlying issue for a range of questions most prominently embodied by M.~Kac's \emph{Can One Hear the Shape of a Drum?}~\cite{Kac66}.
Here we take up this issue and investigate how much of the geometric structure of a space is encoded by an associated Markov semigroup or, equivalently, by a Markov process. To be more precise, let us briefly introduce the setting. For~$i=1,2$, let~$(X^i,\mcX^i,\T^i,\mu^i)$ be a locally compact Polish $\sigma$-finite Radon-measure space and~$\ttseq{T^i_t}_{t\geq 0}$ be a sub-Markovian semigroup on~$L^2(\mu^i)$. They are \emph{intertwined} by $U\colon L^2(\mu^1)\to L^2(\mu^2)$ if
\begin{align*}
U\circ T^1_t=T^2_t\circ U\comma \qquad t\geq 0\fstop
\end{align*}
In~\cite{LenSchWir18} D.~Lenz, M.~Schmidt, and the second named author investigated the properties of intertwined sub-Markovian semigroups in the case when the intertwining operator~$U$ is additionally an \emph{order isomorphism}, that is, an invertible linear operator satisfying $Uf\geq 0$ if and only if~$f\geq 0$.
Suppose that~$\ttseq{T^i_t}_{t\geq 0}$ is the semigroup of a quasi-regular symmetric Dirichlet form~$\ttonde{E^i,\dom{E^i}}$ on~$L^2(\mu^i)$, and thus it is associated with a symmetric Markov process~$\mbfM^i$ with state space~$X^i$.
It is then the main result of~\cite{LenSchWir18} that~$\mbfM^i$ \mbox{(quasi-)}determines the topology~$\T^i$, in the following sense.

\begin{theorem}\label{t:IntroMain}
If the \purple{irreducible} semigroups~$\ttseq{T^i_t}_{t\geq 0}$ are intertwined by an order isomorphism, then the corresponding Dirichlet spaces~$(X^i,\T^i,\mu^i, E^i)$ are quasi-homeo\-morphic.
\end{theorem}
Here, by a \emph{quasi-homeomorphism}\footnote{In the rest of the paper, by a \emph{quasi-homeomorphism} we shall mean a map satisfying the more restrictive conditions in~\cite[Dfn.~3.1]{CheMaRoe94}.}~$j\colon (X^1,\T^1,\mu^1, E^1)\to (X^2,\T^2,\mu^2, E^2)$ we mean a map such that for every~$\eps>0$ there exist closed sets~$F^i\subset X^i$ with $\Cap(X^i\setminus F^i)<\eps$ and~$j$ restricts to a homeomorphism~$j\restr_{F^1}\colon F^1\to F^2$, cf.~\cite[\S{A}.4, p.~429]{FukOshTak11}.

\noindent Importantly, this completely characterizes how a Markov process~$\mbfM$ identifies the topological properties of its state space~$(X,\T,\mu)$, and constitutes a sharp affirmative result in the spirit of M.~Kac' isospectrality question.
For example, ---~as previously established by W.~Arendt, M.~Biegert, and A.~F.~M.~ter Elst in~\cite{AreBietEl12}~--- a Brownian motion on a complete Riemannian manifold~$(M,g)$ characterizes both the differential and the Riemannian structure of~$(M,g)$.

The arguments in~\cite{LenSchWir18} rely however on the additional technical assumption that the semigroups~$\ttseq{T^i_t}_{t\geq 0}$ (or the Markov processes~$\mbfM^i$) be \emph{irreducible}.
In terms of the semigroups this is to say that
\begin{align*}
T_t\car_A = \car_A T_t \text{ for some $t>0$ } \iff  \text{ either } A \text{ or } A^\complement \text{ is } \mu\text{-negligible} \fstop
\end{align*}

Our main contribution in this work is that to extend the results in~\cite{LenSchWir18} by removing the irreducibility assumption, thus proving Theorem~\ref{t:IntroMain} in its full generality.
To this end, we employ the machinery of \emph{ergodic decompositions} of Dirichlet spaces developed by the first named author in~\cite{LzDS20}.
An ergodic decomposition of a Dirichlet space~$\ttonde{E^i,\dom{E^i}}$ is a collection of irreducible Dirichlet spaces~$\ttonde{E_\zeta,\dom{E_\zeta}}$ indexed by a measure space~$(Z,\mcZ,\nu)$ decomposing $\ttonde{E,\dom{E}}$ as the direct integral
\begin{align*}
E= \dint{Z} E_\zeta \diff\nu(\zeta) \fstop
\end{align*}

Extending the uniqueness results in~\cite{LzDS20}, we show that ergodic decompositions are unique up to unique automorphisms of their indexing spaces~$Z$, Theorem~\ref{t:Uniqueness}.
Profiting this fact, we prove as our main result, Theorem~\ref{t:Main}, that two Dirichlet spaces $\ttonde{E^i,\dom{E^i}}$ are intertwined by an order isomorphism~$U$ if and only if their ergodic decompositions~$\zeta\mapsto\ttonde{E^i_\zeta,\dom{E^i_\zeta}}$ are intertwined by an order isomorphism~$U_\zeta$, in which case~$U$ is decomposable over the ergodic decompositions and represented as a direct integral of the operators~$U_\zeta$.

\section{Direct Integrals}\label{s:DirInt}

In this section we review some of the basic theory of direct integrals of Hilbert spaces and quadratic forms for later use. More detailed accounts of the material discussed here can be found in \cite{LzDS20,Dix81}.
We refer to~\cite{Fre00} for all measure-theoretical statements.

\subsection{Measure spaces}
To formulate the ergodic decomposition of Dirichlet forms in full generality, some rather technical aspects of measure theory are required. For the reader's convenience, let us recall the relevant definitions, starting with some classes of measurable spaces that we will use.
\begin{definition}[Measurable spaces]\label{d:Standard} A measurable space~$(X,\mcX)$ is
\begin{itemize}
\item \emph{separable} if~$\mcX$ contains all singletons in~$X$;

\item \emph{countably separated} if there exists a countable family of sets in~$\mcX$ separating points in~$X$;

\item \emph{countably generated} if there exists a countable family of sets in~$\mcX$ generating~$\mcX$ as a $\sigma$-algebra;

\item a \emph{standard Borel space} if there exists a Polish topology~$\T$ on~$X$ so that~$\mcX$ coincides with the Borel $\sigma$-algebra induced by~$\T$.
\end{itemize}

A $\sigma$-finite measure space~$(X,\mcX,\mu)$ is \emph{standard} if there exists a $\mu$-conegligible set $X_0\in\mcX$ such that~$X_0$ is a standard Borel space when regarded as a measurable subspace of~$(X,\mcX)$. We denote by~$(X,\mcX^\mu,\hat\mu)$ the (Carath\'eodory) completion of~$(X,\mcX,\mu)$.
\purple{A Hausdorff topological measure space~$(X,\T,\mcX,\mu)$ is a \emph{Radon-measure space} if~$\mcX^\mu$ coincides with the $\mu$-completion of the Borel $\sigma$-algebra and~$\hat\mu$ is a Radon measure, i.e.\ it is locally finite and inner regular with respect to the compact sets.}

A $[-\infty,\infty]$-valued function is called \emph{$\mu$-measurable} if it is measurable w.r.t.~$\mcX^\mu$. For measures~$\mu_1$,~$\mu_2$ we write~$\mu_1\sim \mu_2$ to indicate that~$\mu_1$ and~$\mu_2$ are equivalent, i.e.\ mutually absolutely continuous.
\end{definition}

Let us now recall the notion of a perfect measure space together with some of its properties.
A measure space~$(X,\mcX,\mu)$ is \emph{perfect} (e.g.~\cite[342K]{Fre00}) if whenever~$f\colon X\to\R$ is measurable and~$A\in\mcX$ with~$\mu A>0$, then there exists a compact set~$K\subset f(A)$ with~$f_\pfwd\mu(K)>0$.
Here, and everywhere in the following, $f_\pfwd\mu$~is the push-forward measure of~$\mu$ via~$f$.

The following  properties of perfect measure spaces are easily verified, cf.~\cite[342X]{Fre00}. For a proof of~\ref{i:l:Perfect:5} see~\cite[416W(a)]{Fre00}.
\begin{lemma}\label{l:Perfect}
Let~$(X,\mcX,\mu)$ be any measure space. Then,
\begin{enumerate}[$(i)$]
\item $(X,\mcX,\mu)$ is perfect if and only if so is its completion;
\item if~$(X,\mcX,\mu)$ is perfect and~$\mcX_0$ is a $\sigma$-subalgebra of~$\mcX$, then the restricted measure space~\mbox{$(X,\mcX_0,\mu\mrestr{\mcX_0})$} is perfect;
\item if~$(X,\mcX,\mu)$ is perfect and~$\lambda$ has a density w.r.t.~$\mu$, then~$(X,\mcX,\lambda)$ is perfect;
\item if~$(Z,\mcZ)$ is a measurable space,~$f\colon X\to Z$ is $\mcX/\mcZ$-measurable and~$(X,\mcX,\mu)$ is perfect, then $(Z,\mcZ,f_\pfwd\mu)$ is perfect;
\item\label{i:l:Perfect:5} every Radon-measure space is perfect.
\end{enumerate}
\end{lemma}

\begin{definition}\label{d:Bianchini}
Let~$(X,\mcX,\mu)$ be a $\sigma$-finite measure space, and~$\mcX^*\subset \mcX$ be a countably generated $\sigma$-subalgebra.
We say that:
\begin{itemize}
\item $\mcX$ is $\mu$-\emph{essentially countably generated by~$\mcX^*$} if for each~$A\in\mcX$ there is~$A^*\in\mcX^*$ with~$\mu(A\triangle A^*)=0$;
\item $\mcX$ is $\mu$-\emph{essentially countably generated} if it is so by some~$\mcX^*$ as above.
\end{itemize}
\end{definition}

\begin{definition}[Almost isomorphisms]\label{d:MeasureIso1}
For~$i\in\{1,2\}$ let~$(X^i,\mcX^i)$ be a measurable space, and~$\mcN^i$ be a $\sigma$-ideal of~$\mcX^i\cup \mcN^i\eqdef \set{A\cup N : A\in\mcX^i, N\in\mcN^i}$. We say that the triples~$(X^i,\mcX^i,\mcN^i)$ are \emph{almost isomorphic} if there exist sets~$N^i\in \mcN^i$ so that the spaces~$X^i\setminus N^i$ are strictly isomorphic when endowed with the restriction to~$X^i\setminus N^i$ of the $\sigma$-algebra~$\mcX^i\cup \mcN^i$ and of the $\sigma$-ideal~$\mcN^i$.

For~$i\in\{1,2\}$ let~$(X^i,\mcX^i,\mu^i)$ be measure spaces. We say that they are \emph{almost isomorphic} if there exist $\mu^i$-negligible sets~$N^i\subset X^i$ so that the spaces~$X^i\setminus N^i$ are strictly isomorphic when endowed with the restriction to~$X^i\setminus N^i$ of the completed $\sigma$-algebra~$(\mcX^i)^{\mu^i}$, of the completed measure~$\hat\mu^i$, and of the $\sigma$-ideal of~$\mu^i$-negligible sets~$\mcN^{\mu^i}$.
\end{definition}

\subsection{Quadratic forms}
Every Hilbert space is assumed to be \emph{separable} and a real vector space.

Let~$(H,\norm{\emparg})$ be a Hilbert space.
\purple{For a bounded operator~$B\colon H\to H$ we write~$\norm{B}_{\op, H}$ for the operator norm of~$B$; the subscript~$H$ is omitted whenever apparent from context.}
A \emph{quadratic form}~$Q$ on~$H$ is a map $Q\colon H\to [0,\infty]$ such that
\begin{equation*}
Q(u+v)+Q(u-v)=2Q(u)+2Q(v)\comma\qquad u,v\in H.
\end{equation*}
We always assume that the domain $\dom{Q}=\{u\in H\colon Q(u)<\infty\}$ is dense in~$H$ and write $(Q,\dom{Q})$ for~$Q$ when we want to make the domain explicit.

For every quadratic form~$Q$ there exists a unique symmetric bilinear form whose diagonal coincides with~$Q$. This bilinear form will also be denoted by~$Q$, and we will use the term quadratic form for both objects interchangeably.

Additionally, for every~$\alpha>0$ we set
\begin{equation*}
Q_\alpha(u)\eqdef Q(u)+\alpha\norm{u}^2\comma \qquad u\in H \fstop
\end{equation*}
For~$\alpha>0$, we let~$\dom{Q}_\alpha$ be the completion of~$\dom{Q}$, endowed with the Hilbert norm~$Q_\alpha^{1/2}$.

To every closed quadratic form~$(Q,\dom{Q})$ one can associate a unique non-negative self-adjoint operator~$-L$ such that $\dom{\sqrt{-L}}=\dom{Q}$ and $Q(u,v) = \scalar{-Lu}{v}$ for all~$u,v\in\dom{L}$. We denote the associated strongly continuous contraction semigroup by~$T_t\eqdef e^{tL}$,~$t>0$.%, and the associated strongly continuous contraction resolvent by~$G_\alpha\eqdef (\alpha-L)^{-1}$,~$\alpha> 0$.
%By the Hille--Yosida Theorem, e.g.~\cite[p.~27]{MaRoe92},
%\begin{subequations}\label{eq:Hille--Yosida}
%\begin{align}
%\label{eq:Hille--YosidaG}
%Q_\alpha(G_\alpha u,v)=&\ \scalar{u}{v}_H\comma\qquad u\in H\comma v\in\dom{Q}\comma 
%\\
%\label{eq:Hille--YosidaT}
%T_tu=&\lim_{\alpha\rar \infty} e^{t\alpha (\alpha G_\alpha-1)}u,\qquad u\in H\comma
%\\
%\label{eq:Hille--YosidaE}
%Q(u,v)=&\ \lim_{\beta\rar\infty} \scalar{\beta u-\beta G_\beta u}{v}_H\comma \qquad u,v\in H
%\fstop
%\end{align}
%\end{subequations}

\subsection{Direct integrals of Hilbert spaces}
We recall the main definitions concerning direct integrals of separable Hilbert spaces, referring to~\cite[\S\S{II.1}, II.2]{Dix81} for a systematic treatment.

\begin{definition}[Measurable fields,~{\cite[\S{II.1.3}, Dfn.~1,~p.~164]{Dix81}}]\label{d:DirInt}
\
Let~$(Z,\mcZ,\nu)$ be a $\sigma$-finite measure space,~$\seq{H_\zeta}_{\zeta\in Z}$ be a family of separable Hilbert spaces, and~$F$ be the linear space~$F\eqdef \prod_{\zeta\in Z} H_\zeta$. We say that~$\zeta\mapsto H_\zeta$ is a \emph{$\nu$-measurable field of Hilbert spaces} (\emph{with underlying space~$S$}) if there exists a linear subspace~$S$ of~$F$ with
\begin{enumerate}[$(a)$]
\item\label{i:d:DirInt1} for every~$u\in S$, the function~$\zeta\mapsto \norm{u_\zeta}_{\zeta}$ is $\nu$-measurable;
\item\label{i:d:DirInt2} if~$v\in F$ is such that $\zeta\mapsto \scalar{u_\zeta}{v_\zeta}_{\zeta}$ is $\nu$-measurable for every~$u\in S$, then~$v\in S$;
\item\label{i:d:DirInt3} there exists a sequence~$\seq{u_n}_n\subset S$ such that~$\seq{u_{n,\zeta}}_n$ is a total sequence\footnote{A sequence in a Banach space~$B$ is called \emph{total} if the strong closure of its linear span coincides with~$B$.} in~$H_\zeta$ for every~$\zeta\in Z$.
\end{enumerate}

Any such $S$ is called a \emph{space of $\nu$-measurable vector fields}. Any sequence in~$S$ possessing property~\iref{i:d:DirInt3} is called a \emph{fundamental} sequence.
\end{definition}

\begin{proposition}[{\cite[\S{II.1.4, Prop.~4, p.~167}]{Dix81}}]\label{p:Dix} Let~$\mcS$ be a family of functions satisfying Definition~\ref{d:DirInt} \iref{i:d:DirInt1} and~\iref{i:d:DirInt3}. Then, there exists exactly one space of $\nu$-measurable vector fields~$S$ so that~$\mcS\subset S$.
\end{proposition}

\begin{definition}[Direct integrals, {\cite[\S{II.1.5}, Prop.~5,~p.~169]{Dix81}}]
A $\nu$-measur\-able vector field~$u$ is called ($\nu$-)\emph{square-integrable} if
\begin{align}\label{eq:NormH}
\norm{u}\eqdef \tonde{\int_Z \norm{u_\zeta}_{\zeta}^2 \, \diff\nu(\zeta)}^{1/2}<\infty\fstop
\end{align}

Two square-integrable vector fields~$u$,~$v$ are called ($\nu$-)\emph{equivalent} if~$\norm{u-v}=0$. The space~$H$ of equivalence classes of square-integrable vector fields, endowed with the non-relabeled quotient norm~$\norm{\emparg}$, is a Hilbert space~\cite[\S{II.1.5}, Prop.~5(i), p.~169]{Dix81}, called the \emph{direct integral of~$\zeta\mapsto H_\zeta$} (\emph{with underlying space~$S$}) and denoted by
\begin{align}\label{eq:DirInt}
H=\dint[S]{Z} H_\zeta \diff\nu(\zeta) \fstop
\end{align}
The superscript~`$S$' is omitted whenever~$S$ is clear from context.
\end{definition}

In the following, it will occasionally be necessary to distinguish an element~$u$ of~$H$ from one of its representatives modulo $\nu$-equivalence, say~$\rep u$ in~$S$. In this case, we shall write~$u=\class[H]{\rep u}$.

\begin{lemma}[{\cite[Lem.~2.6]{LzDS20}}] If~$(Z,\mcZ)$ is $\sigma$-finite countably generated, then~$H$ in~\eqref{eq:DirInt} is separable.
\end{lemma}

We now turn to measurable fields of bounded operators.

\begin{definition}[Measurable fields of bounded operators, decomposable operators]\label{d:BoundedOps}
Let~$H$ be as in~\eqref{eq:DirInt}.
A field of bounded operators~$\zeta\mapsto B_\zeta\in\mcB(H_\zeta)$ is called \emph{$\nu$-measurable} (\emph{with underlying space~$S$}) if~$\zeta\mapsto B_\zeta u_\zeta\in H_\zeta$ is a $\nu$-measurable vector field for every $\nu$-measurable vector field~$u$.
\purple{Set $\norm{B_\zeta}_{\op, \zeta}\eqdef \norm{B_\zeta}_{\op, H_\zeta}$.}
A $\nu$-measurable vector field of bounded operators is called \emph{$\nu$-essentially bounded} if~$\nu$-$\esssup_{\zeta\in Z}\norm{B_\zeta}_{\op,\zeta}<\infty$. A bounded operator~$B\in\mcB(H)$ is called \emph{decomposable} if~$Bu$ is represented by a $\nu$-essentially bounded $\nu$-measurable field of bounded operators~$\zeta\mapsto B_\zeta$, in which case we write
\begin{align*}
B=\dint{Z} B_\zeta \diff\nu(\zeta) \fstop
\end{align*}
\end{definition}

\subsection{Direct integrals of quadratic forms}\label{ss:DirIntQ}
We briefly recall all the relevant notions concerning direct integrals of quadratic forms according to~\cite{LzDS20}.

\begin{definition}[Direct integral of quadratic forms]\label{d:assQ}
Let~$(Z,\mcZ,\nu)$ be a $\sigma$-finite countably generated measure space. For~$\zeta\in Z$ let~$(Q_\zeta,D_\zeta)$ be a closable (densely defined) quadratic form on a Hilbert space~$H_\zeta$. We say that~$\zeta\mapsto (Q_\zeta,D_\zeta)$ is a \emph{$\nu$-measurable field of quadratic forms on~$Z$} if
\begin{enumerate}[$(a)$]
\item\label{i:assQ:a} $\zeta\mapsto H_\zeta$ is a $\nu$-measurable field of Hilbert spaces on~$Z$ with underlying space~$S_H$;
\item\label{i:assQ:b} $\zeta\mapsto \dom{Q_\zeta}_1$ is a $\nu$-measurable field of Hilbert spaces on~$Z$ with underlying space~$S_Q$;
\item\label{i:assQ:c} $S_Q$ is a linear subspace of~$S_H$ under the canonical identification of~$\dom{Q_\zeta}$ with a subspace of~$H_\zeta$.
\end{enumerate}

We \purple{further} denote by
\begin{align*}
Q=\dint[S_Q]{Z} Q_\zeta \diff\nu(\zeta)
\end{align*}
the \emph{direct integral} of~$\zeta\mapsto (Q_\zeta,\dom{Q_\zeta})$, i.e. the quadratic form defined on~$H$ as in~\eqref{eq:DirInt} given by
\begin{equation}\label{eq:DIntQF}
\begin{aligned}
\dom{Q}\eqdef& \set{\class[H]{\rep u} : \rep u \in S_Q , \int_Z Q_{\zeta,1}(\rep u_\zeta) \diff\nu(\zeta) <\infty} \comma
\\
Q(u,v)\eqdef& \int_Z Q_\zeta(u_\zeta,v_\zeta) \diff\nu(\zeta)\comma \qquad u,v\in\dom{Q}\fstop
\end{aligned}
\end{equation}
\end{definition}

\begin{remark}[Separability]\label{r:Separability} It is implicit in our definition of $\nu$-measurable field of Hilbert spaces that~$H_\zeta$ is \emph{separable} for every~$\zeta\in Z$. As a consequence, $\dom{Q_\zeta}_1$ is $(Q_\zeta)^{1/2}_1$-separable.
\end{remark}

\begin{proposition}[{\cite[Prop.~2.13]{LzDS20}}]\label{p:DirInt}
Let~$(Q,\dom Q)$ be a direct integral of quadratic forms. Then,
\begin{enumerate}[$(i)$]
\item\label{i:p:DirInt1} $\ttonde{Q,\dom Q}$ is a densely defined closed quadratic form on~$H$;
\item\label{i:p:DirInt2} $\zeta\mapsto T_{\zeta,t}$ is $\nu$-measurable fields of bounded operators for every~$t>0$;
\item\label{i:p:DirInt3} the semigroup associated with $Q$ is given by
\begin{equation}\label{eq:p:DirInt:0}
T_t \eqdef \dint[S_H]{Z} T_{\zeta,t} \diff\nu(\zeta) \comma \quad t>0\fstop
\end{equation}
\end{enumerate}
\end{proposition}

Under the assumptions of Proposition~\ref{p:DirInt}, assertion~\ref{i:p:DirInt1} of the same Proposition implies that the space of $\nu$-measurable vector fields~$S_H$ is uniquely determined by~$S_Q$ as a consequence of Proposition~\ref{p:Dix}. Thus, everywhere in the following when referring to a direct integral of quadratic forms we shall ---~with abuse of notation~--- write $S$ in place of both~$S_H$ and~$S_Q$.

\subsection{Direct-integral representation of \texorpdfstring{$L^2$}{L2}-spaces}\label{ss:DirIntE}
In order to introduce direct-integral representations of Dirichlet forms, we need to construct direct integrals of concrete Hilbert spaces in such a way to additionally preserve the Riesz structure of Lebesgue spaces implicitly used to phrase the sub-Markovianity property~\eqref{eq:subMarkov}. To this end, we shall need the concept of a disintegration of measures.

\paragraph{Disintegrations} Let~$(X,\mcX,\mu)$ and~$(Z,\mcZ,\nu)$ be (non-trivial) measure spaces. A map~$s\colon (X,\mcX)\rar (Z,\mcZ)$ is \emph{inverse-measure-preserving} if~$s_\pfwd \mu\eqdef \mu \circ s^{-1}=\nu$.

\begin{definition}[Disintegrations~{\cite[452E]{Fre00}}]\label{d:Disint}
A \emph{pseudo-disintegration of~$\mu$ over~$\nu$} is any family of non-zero measures~$\seq{\mu_\zeta}_{\zeta\in Z}$ on~$(X,\mcX)$ so that~$\zeta\mapsto \mu_\zeta A$ is $\nu$-measurable and
\begin{align*}
\mu A=\int_Z\mu_\zeta A \diff\nu(\zeta)\comma \qquad A\in\mcX\fstop
\end{align*}
A pseudo-disintegration is
\begin{itemize}
\item \emph{separated} if there exists a family of pairwise disjoint sets~$\set{A_\zeta}_{\zeta\in Z}\subset \mcX^\mu$ so that~$A_\zeta$ is $\mu_\zeta$-conegligible for $\nu$-a.e.~$\zeta\in Z$, henceforth called a \emph{separating family} for~$\seq{\mu_\zeta}_{\zeta\in Z}$;

\item \emph{$s$-separated} if it is separated and there exists a $\mcX/\mcZ$-measurable map~$s\colon X\rar Z$ so that~$\set{s^{-1}(\zeta)}_{\zeta\in Z}$ is a separating family;

\item \emph{consistent with~$s$} if
\begin{align}\label{eq:Disint:0}
\mu \ttonde{A\cap s^{-1}(B)}=\int_B\mu_\zeta A \diff\nu(\zeta)\comma \qquad A\in\mcX\comma B\in \mcZ
\end{align}

\item \emph{strongly consistent with~$s$} if it is consistent  with~$s$ and $s^{-1}(\zeta)$ is~$\mu_\zeta$-coneglig-ible for $\nu$-a.e.~$\zeta\in Z$.
\end{itemize}

A \emph{disintegration of~$\mu$ over~$\nu$} is a pseudo-disintegration additionally so that~$\mu_\zeta$ is a sub-probability measure for every~$\zeta\in Z$.
A disintegration is
\begin{itemize}
\item $\nu$-\emph{essentially unique} if the measures~$\mu_\zeta$ are uniquely determined for $\nu$-a.e.~$\zeta\in Z$.
\end{itemize}
\end{definition}

If~$\seq{\mu_\zeta}_{\zeta\in Z}$ is a pseudo-disintegration of~$\mu$ over~$\nu$, then
\begin{align}\label{eq:DisintF}
\int_X g\diff\mu=\int_Z\int_X g(x) \diff\mu_\zeta(x) \diff\nu(\zeta)
\end{align}
whenever the left-hand side makes sense,~\cite[452F]{Fre00}. We note that a disintegration~$\seq{\mu_\zeta}_{\zeta\in Z}$ of~$\mu$ over~$\nu$ strongly consistent with a map~$s$ is automatically separated, with separating family~$\ttset{s^{-1}(\zeta)}_{\zeta\in Z}$.

\paragraph{Direct integrals and disintegrations} Let~$(X,\mcX,\mu)$ be $\sigma$-finite standard,~$(Z,\mcZ,\nu)$ be $\sigma$-finite countably generated, and~$\seq{\mu_\zeta}_{\zeta\in Z}$ be a pseudo-disintegration of~$\mu$ over~$\nu$. Denote by
\begin{itemize}
\item $\mcL^0(\mu)$ the space of $\mu$-measurable real-valued functions (\emph{not}: $\mu$-classes) on~$X$;
\item $\mcL^\infty(\mu)$ the space of uniformly bounded (\emph{not}: $\mu$-essentially uniformly bound\-ed) functions in~$\mcL^0(\mu)$;
\item $\mcL^p(\mu)$ the space of $p$-integrable functions in~$\mcL^0(\mu)$.
\end{itemize}

For a family~$\mcA\subset\mcL^0(\mu)$, let~$\class[\mu]{\mcA}$ denote the family of the corresponding $\mu$-classes.

\smallskip

Let now~$F\eqdef \prod_{\zeta\in Z} L^2(\mu_\zeta)$. The \emph{diagonal embedding} of~$\mcL^2(\mu)$ into~$F$, regarded up to $\mu_\zeta$-classes, is the map~$\delta\colon f\mapsto \ttseq{\zeta\mapsto \delta(f)_\zeta}$, where
\begin{align}\label{eq:Delta}
\delta(f)_\zeta\eqdef \begin{cases} \class[\mu_\zeta]{f} & \text{if~} f\in \mcL^2(\mu_\zeta), \\ \zero_{L^2(\mu_\zeta)} & \text{otherwise}.\end{cases}
\end{align}
Since $f\in \mcL^2(\mu)$, we have $\delta(f)_\zeta=\class[\mu_\zeta]{f}$ for~$\nu$-a.e.~$\zeta\in Z$ by~\eqref{eq:DisintF}, and~$\delta$ therefore is well-defined as linear morphism mapping~$\mu$-classes to $H$-classes, see Proposition~\ref{p:DIntL^2} below.
Now, assume that
\begin{equation}\label{eq:AssA}
\text{$\mcA$ is a linear subspace of~$\mcL^2(\mu)$, and~$\class[\mu]{\mcA}$ is dense in~$L^2(\mu)$.}
\end{equation}
Since~$\class[\mu]{\mcA}$ is dense in $L^2(\mu)$ and the latter is separable, then there exists a countable family~$\mcU\subset \mcA$ so that~$\class[\mu_\zeta]{\mcU}$ is total in~$L^2(\mu_\zeta)$ for $\nu$-a.e.~$\zeta\in Z$.
Thus for every~$\mcA$ as in~\eqref{eq:AssA} there exists a unique space of $\nu$-measurable vector fields~$S=S_\mcA$ containing~$\delta(\mcA)$, generated by~$\delta(\mcA)$ in the sense of Proposition~\ref{p:Dix}. We denote by~$H$ the corresponding direct integral of Hilbert spaces
\begin{align}\label{eq:H}
H\eqdef \dint[S]{Z} L^2(\mu_\zeta) \diff\nu(\zeta)\fstop
\end{align}

\begin{remark}[{cf.~\cite[\S7.2, p.~84]{HayMirYve91}}]\label{r:Order}
The direct integral $H$ constructed in~\eqref{eq:H} is a \emph{Banach lattice} (e.g.~\cite[354A(b)]{Fre00}) for the order
\begin{align*}
h\geq  \zero_H \qquad \iff \qquad h_\zeta\geq \zero_{L^2(\mu_\zeta)} \quad \forallae{\nu} \zeta\in Z \fstop
\end{align*}
\end{remark}

\begin{proposition}[{\cite[Prop.~2.25]{LzDS20}}]\label{p:DIntL^2}
Let~$(X,\mcX,\mu)$ be $\sigma$-finite standard, $(Z,\mcZ,\nu)$ be $\sigma$-finite countably generated, and~$\seq{\mu_\zeta}_{\zeta\in Z}$ be a pseudo-disintegration of~$\mu$ over~$\nu$. Then, the morphism
\begin{align}\label{eq:DirIntL^2:0}
\iota\colon L^2(\mu) \longrar H\eqdef \dint[S]{Z} L^2(\mu_\zeta) \diff\nu(\zeta)\comma \qquad \class[\mu]{f} \longmapsto \class[H]{\delta(f)}
\end{align}
\begin{enumerate}[$(i)$]
\item\label{i:p:DIntL^2:1} is well-defined, linear, and an isometry of Hilbert spaces, additionally unitary if~$\seq{\mu_\zeta}_{\zeta\in Z}$ is separated;

\item\label{i:p:DIntL^2:2} is a Riesz homomorphism (e.g.~\cite[351H]{Fre00}). In particular,
\begin{itemize}
\item for each~$f\in\mcL^2(\mu)$, it holds that~$(\iota\class[\mu]{f})_\zeta \geq \zero_{L^2(\mu_\zeta)}$ for~$\nu$-a.e.~$\zeta\in Z$ if and only if~$f\geq 0$ $\mu$-a.e.;

\item for each~$f,g\in\mcL^2(\mu)$, it holds that~$(\iota\class[\mu]{f \wedge g})_\zeta= (\iota\class[\mu]{f})_\zeta \wedge (\iota\class[\mu]{g})_\zeta$ for~$\nu$-a.e.~$\zeta\in Z$.
\end{itemize}
\end{enumerate}
\end{proposition}

\subsection{Dirichlet forms}
We recall a standard setting for the theory of Dirichlet forms, following~\cite{MaRoe92}.

\begin{Assumption}\label{ass:Main} The quadruple $(X,\T,\mcX,\mu)$ is so that~$(X,\T)$ is a metrizable Luzin space with Borel $\sigma$-algebra~$\mcX$ and~$\hat\mu$ is a Radon measure on~$(X,\T,\mcX^\mu)$ with full support.
\end{Assumption}

By~\cite[415D(iii), 424G]{Fre00} any space~$(X,\mcX,\mu)$ satisfying Assumption~\ref{ass:Main} is $\sigma$-finite standard. The \emph{support} of a ($\mu$-)measurable function~$f\colon X\rar \R$ (possibly defined only on a $\mu$-conegligible set) is defined as the support of the measure $\abs{f}\cdot \mu$. Every such~$f$ has a support, independent of the $\mu$-representative of~$f$, cf.~\cite[p.~148]{MaRoe92}.

A closed positive semi-definite quadratic form~$(Q,\dom{Q})$ on~$L^2(\mu)$ is a (\emph{symmetric}) \emph{Dirichlet form} if
\begin{align}\label{eq:subMarkov}
f\in\dom{Q} \implies f^+ \wedge \car \in \dom{Q} \text{~and~} Q(f^+\wedge \car)\leq Q(f) \fstop
\end{align}
We shall denote Dirichlet forms by~$(E,\dom{E})$. A Dirichlet form~$(E,\dom{E})$ is \emph{regular} if~$(X,\T)$ is (additionally) locally compact and $\dom{E}\cap \Cz(X)$ is both $E^{1/2}_1$-dense in~$\dom{E}$ and uniformly dense in~$\Cz(X)$.

Finally, we are interested in the notion of invariant sets of a Dirichlet form~$(E,\dom{E})$ on~$L^2(\mu)$.
We say that~$A\subset X$ is $E$-\emph{invariant} if it is $\mu$-measurable and any of the following equivalent %\footnote{See~\cite[Lem.~1.6.1, p.~53]{FukOshTak11}, the proof of which adapts \emph{verbatim} to our more general setting.}
conditions holds.
\begin{itemize}
\item $T_t(\car_A f)=\car_A T_t f$ $\mu$-a.e. for any~$f\in L^2(\mu)$ and~$t>0$;
\item $T_t(\car_A f)=0$ $\mu$-a.e. on $A^\complement$ for any~$f\in L^2(\mu)$ and~$t>0$;
%\item $G_\alpha(\car_A f)=0$ $\mu$-a.e. on $A^\complement$ for any~$f\in L^2(\mu)$ and~$\alpha>0$;
\item $\car_A f\in \dom{E}$ for any~$f\in \dom{E}$ and
\begin{align}\label{eq:d:Invariant1}
E(f,g)=E(\car_A f,\car_A g)+E(\car_{A^\complement}f,\car_{A^\complement} g)\comma \qquad f,g\in \dom{E} \fstop
\end{align}
\end{itemize}
The form~$(E,\dom{E})$ is \emph{irreducible} if every invariant set is either negligible or conegligible.
%\end{definition}

\subsection{Direct integrals of Dirichlet forms} Let~$(X,\mcX,\mu)$ be $\sigma$-finite standard,~$(Z,\mcZ,\nu)$ be $\sigma$-finite countably generated, and~$\seq{\mu_\zeta}_{\zeta\in Z}$ be a pseudo-disintegration of~$\mu$ over~$\nu$. Further let~$\zeta\mapsto(Q_\zeta, \dom{Q_\zeta})$ be a $\nu$-measurable field of quadratic forms, each densely defined in~$L^2(\mu_\zeta)$ with separable domain, and denote by~$(Q,\dom{Q})$ their direct integral in the sense of Definition~\ref{d:assQ}.

\begin{definition}[{\cite[Dfn.~2.26]{LzDS20}}]\label{d:Compat}
We say that~$(Q,\dom{Q})$ is \emph{compatible} with the pseudo-disintegration $\seq{\mu_\zeta}_{\zeta\in Z}$ if the space~$S_Q$ underlying~$\zeta\mapsto \dom{Q_\zeta}_1$ is of the form~$S_\mcA$ for some~$\mcA$ as in~\eqref{eq:AssA} and additionally satisfying~$\mcA\subset \dom{Q}$.
\end{definition}

\begin{proposition}[{\cite[Prop.~2.29]{LzDS20}}]\label{p:DirIntE}
Let~$(X,\mcX,\mu)$ be $\sigma$-finite standard,~$(Z,\mcZ,\nu)$ be $\sigma$-finite countably generated, and~$\seq{\mu_\zeta}_{\zeta\in Z}$ be a separated pseudo-disintegration of~$\mu$ over~$\nu$. Further let~$(E,\dom{E})$ be a direct integral of quadratic forms~$\zeta\mapsto(E_\zeta, \dom{E_\zeta})$ compatible with~$\seq{\mu_\zeta}_{\zeta\in Z}$. Then, $(E,\dom{E})$ is a Dirichlet form on~$L^2(\mu)$ if and only if~$(E_\zeta, \dom{E_\zeta})$ is so on~$L^2(\mu_\zeta)$ for $\nu$-a.e.~$\zeta\in Z$.
\end{proposition}

Proposition~\ref{p:DirIntE} motivates the following definition.

\begin{definition}\label{d:DirIntE}
A quadratic form~$(E,\dom{E})$ on~$L^2(\mu)$ is a \emph{direct integral of Dirichlet forms}~$\zeta\mapsto (E_\zeta,\dom{E_\zeta})$ if it is a direct integral \purple{of quadratic forms}~$\zeta\mapsto (E_\zeta,\dom{E_\zeta})$, \purple{for each~$\zeta$ the form~$(E_\zeta,\dom{E_\zeta})$ is a Dirichlet form on~$L^2(\mu_\zeta)$, and the direct integral is} additionally compatible with the \emph{separated} pseudo-disintegration~$\seq{\mu_\zeta}_\zeta$ in the sense of Definition~\ref{d:Compat}.
\end{definition}

We refer to~\cite{LzDS20} for further comments on direct integrals of Dirichlet forms and related notions.

\section{Order-isomorphisms and intertwining operators}

In this section we first recall the notion of order isomorphisms between $L^2$-spaces together with their main representation theorem, which shows that they are weighted composition operators. In the second part we study intertwining operators and give a first result connecting intertwining operators and direct integral decompositions (Proposition \ref{p:IntertwiningDirInt}).

\begin{definition}[Order-preserving operators, order isomorphisms]
For $i=1,2$ let~$(X^i,\mcX^i,\mu^i)$ be $\sigma$-finite countably generated \purple{and countably separated}.
A linear operator $U\colon L^2(\mu^1)\rar L^2(\mu^2)$ is \emph{order-preserving} (also: \emph{positiv\-ity-preserving}) if~$f\geq 0$ implies~$Uf\geq 0$ for each~$f\in L^2(\mu^1)$.
An \emph{order isomorphism} is an invertible order-preserving linear operator with order-preserving inverse.
\end{definition}

The structure of order isomorphisms between~$L^2$-spaces is characterized by the following Banach--Lamperti-type theorem.

\begin{proposition}[Order isomorphism as weighted composition operator] \label{p:Lamperti} If \sloppy$U\colon L^2(X^1,\mu^1)\rar L^2(X^2,\mu^2)$ is an order isomorphism, then there exist a measurable map~$h\colon X^2 \rar (0,\infty)$ and an $\mcX^2/\mcX^1$-measurable almost isomorphism~$\tau\colon X^2\rar X^1$ such
that
\begin{align*}
Uf=h\cdot (f\circ\tau)\comma \qquad f\in L^2(X^1,\mu^1) \fstop
\end{align*}
The maps~$h$ and~$\tau$ are unique up to equality almost everywhere.
\end{proposition}
\begin{proof}
See~\cite[Thm.~5.1]{Wei84} \purple{for existence in} the case of finite measures~$\mu^i$.
For the general case argue as follows. Since~$(X^i,\mcX^i,\mu^i)$ is $\sigma$-finite standard, there exists a function~$f_i$ satisfying~$f_i>0$ $\mu^i$-a.e., and~$\norm{f_i}_{L^2(\mu^i)}=1$.
Thus, the multiplication operator~$M_i\eqdef M_{f_i^{-1/2}}\colon L^2(\mu^i)\to L^2(f_i\mu^i)$ is an order isomorphism with (order) inverse~$M_i^{-1}=M_{f_i^{1/2}}$, and therefore the operator~$U'\eqdef M_2\circ U\circ M_1^{-1}\colon L^2(f_1\mu^1)\to L^2(f_2\mu^2)$ is an order isomorphism as well.
Applying the assertion for~$U'$ yields maps~$h'$ and~$\tau'$.
Finally, letting~$h\eqdef f_2^{-1/2}\circ \tau \cdot f_1^{1/2} \cdot h'$ and~$\tau\eqdef \tau'$ yields the assertion in the $\sigma$-finite case.

\purple{For uniqueness, one can similarly reduce to the case of finite measure. Then $h=U1$ is determined uniquely up to equality a.e.
Moreover, if $\tau_1,\tau_2\colon X^2\to X^1$ are $\mathcal X^2$/$\mathcal X^1$-measurable maps such that $h\cdot(f\circ\tau_1)=h\cdot(f\circ\tau_2)$ for all $f\in L^2(X^1,\mu^1)$, then $\mu^1(\tau_1^{-1}(B)\cap \tau_2^{-1}(B))=0$ for all $B\in \mathcal X^2$. From~\cite[343F]{Fre00} it follows that $\tau_1=\tau_2$ $\mu^2$-a.e.}
\qed\end{proof}

The maps~$h$ and $\tau$ associated with an order isomorphism~$U$ according to the previous proposition are the main players in the present article. We call~$h$ the \emph{scaling} and~$\tau$ the \emph{transformation} associated with~$U$.

\begin{lemma}[Adjoint of an order isomorphism,~{\cite[Lem.~1.2]{LenSchWir18}}]\label{l:Adjoint} Let
$U\colon L^2(X^1,\mu^1)\rar L^2(X^2,\mu^2)$ be an order isomorphism with associated scaling~$h$ and transformation~$\tau$.
Then $\tau_\pfwd \mu^2$ and $\mu^1$ are mutually absolutely continuous. Moreover, the adjoint
of~$U$ is given by
\begin{align*}
U^*\colon L^{2}(X^2,\mu^2)\rar L^{2}(X^1,\mu^1)\comma \qquad U^* g=\frac{\diff(\tau_\pfwd \mu^2)}{\diff \mu^1}(hg)\circ\tau^{-1}\comma \qquad g\in L^2(\mu^2)\fstop
\end{align*}
\end{lemma}

\begin{definition}[Intertwining operators]
For~$i=1,2$ let~$V^\sym{i}$ be a Banach space and~$(A^\sym{i},\dom{A^\sym{i}})$ be a (possibly unbounded) densely defined operator on~$V^\sym{i}$. An invertible bounded linear operator~$U\colon V^\sym{1}\rar V^\sym{2}$ is said to \emph{intertwine} $A^\sym{1}$ and $A^\sym{2}$ if
\begin{align*}
U\dom{A^\sym{1}}=\dom{A^\sym{2}} \qquad \text{and} \qquad U A^\sym{1} f=A^\sym{2} U f\comma \quad f\in \dom{A^\sym{1}}\fstop
\end{align*}
Two families of operators~$\seq{A^\sym{i}_t}_{t\in T}$,~$i=1,2$, are \emph{intertwined} by~$U$ if~$U$ intertwines~$A^\sym{1}_t$,~$A^\sym{2}_t$ for all~$t\in T$.
\end{definition}

The following is not difficult to show.

\begin{lemma}[{\cite[Prop.~A.1]{KelLenSchWir15}}] Two strongly continuous contraction semigroups on Banach spaces are intertwined by~$U$ if and only if so are their generators.
\end{lemma}

The next result is a variation of~\cite[Cor.~2.5]{LenSchWir18}, replacing the irreducibility assumption on the semigroups with the unitariness of the intertwining operator~$U$.
A proof is analogous, and therefore it is omitted.
\begin{proposition}\label{p:IntertwiningDirichlet}
For~$i=1,2$ let~$(X^i,\mcX^i,\mu^i)$ be $\sigma$-finite standard,~$\seq{T^\sym{i}_t}_{t\geq 0}$ be a sub-Markovian semigroup on~$L^2(\mu^i)$ with corresponding Dirichlet form $(E^i,\dom{E^i})$, and~$U\colon L^2(\mu^1)\rar L^2(\mu^2)$ be a unitary order-isomorphism intertwining~$\seq{T^\sym{1}_t}_{t\geq 0}$ and~$\seq{T^\sym{2}_t}_{t\geq 0}$.
Then,
\begin{align*}
U\dom{E^1}=\dom{E^2} \qquad \text{and} \qquad E^2(U f, U g)=E^1(f,g)\comma \qquad f,g\in\dom{E^1} \fstop
\end{align*}
\end{proposition}

\begin{proposition}\label{p:IntertwiningDirInt}
For~$i=1,2$ let~$(X^i,\mcX^i,\mu^i)$ be $\sigma$-finite standard,~$(Z,\mcZ,\nu)$ be $\sigma$-finite countably generated, and~$\tseq{\mu^i_\zeta}_{\zeta\in Z}$ be a separated pseudo-disintegration of~$\mu^i$ over~$\nu$.
Further let~$\zeta\mapsto U_\zeta$ be a $\nu$-measurable field of bounded operators~$U_\zeta\colon L^2(\mu^1_\zeta)\rar L^2(\mu^2_\zeta)$ in the sense of Definition~\ref{d:BoundedOps}, and set
\begin{align*}
U\eqdef \dint[]{Z} U_\zeta\diff\nu(\zeta) \fstop
\end{align*}

Then,~$U$ is a bounded operator $U\colon L^2(\mu^1)\rar L^2(\mu^2)$ and
\begin{enumerate}[$(i)$]
\item\label{i:p:IntertwiningDirInt:1} $U\colon L^2(\mu^1)\rar L^2(\mu^2)$ is unitary if and only if~$U_\zeta$ is so for $\nu$-a.e.~$\zeta\in Z$, in which case
\begin{align*}
U^{-1}=\dint[]{Z} U_\zeta^{-1}\diff\nu(\zeta)\fstop
\end{align*}

\item\label{i:p:IntertwiningDirInt:2} $U\colon L^2(\mu^1)\rar L^2(\mu^2)$ is order-preserving if and only if~$U_\zeta$ is so for $\nu$-a.e.~$\zeta\in Z$.
\end{enumerate}

For~$i=1,2$ further let~$\zeta\mapsto Q_\zeta^i$ be a $\nu$-measurable field of closed quadratic forms~$\ttonde{Q_\zeta^i,\dom{Q^i_\zeta}}$ on~$L^2(\mu^i_\zeta)$ in the sense of Definition~\ref{d:assQ}, each with semigroup~$\tseq{T^\sym{i}_{\zeta,t}}_{t\geq 0}$, $\zeta\in Z$, and set
\begin{align*}
T^\sym{i}_t\eqdef \dint[]{Z}T^\sym{i}_{\zeta,t}\diff\nu(\zeta) \comma \qquad t>0\fstop
\end{align*}

Then, for each fixed~$t>0$, the operator~$T^\sym{i}_t$ satisfies $T^\sym{i}_t\in\mcB\ttonde{L^2(\mu^i)}$ and
\begin{enumerate}[$(i)$]\setcounter{enumi}{2}
\item\label{i:p:IntertwiningDirInt:3} $U\colon L^2(\mu^1)\rar L^2(\mu^2)$ intertwines~$T^\sym{1}_t$ and $T^\sym{2}_t$ if and only if~$U_\zeta$ intertwines~$T^\sym{1}_{\zeta,t}$ and $T^\sym{2}_{\zeta,t}$ for $\nu$-a.e.~$\zeta\in Z$.
\end{enumerate}
\end{proposition}

\begin{proof} By definition of~$U$, we have
\begin{align*}
U\colon \dint[]{Z} L^2(\mu^1_\zeta) \diff\nu(\zeta)\longrar  \dint[]{Z} L^2(\mu^2_\zeta) \diff\nu(\zeta) \fstop
\end{align*}
Let~$i=1,2$. Since the disintegration~$\ttseq{\mu^i_\zeta}_{\zeta\in Z}$ is separated, by Proposition~\ref{p:DIntL^2} \iref{i:p:DIntL^2:1} there exists a unitary order-isomorphism~$\iota^\sym{i}$ satisfying~\eqref{eq:DirIntL^2:0}.
Since~$\iota^\sym{i}$ is as well an order isomorphism by Proposition~\ref{p:DIntL^2}~\iref{i:p:DIntL^2:2}, in the rest of the proof we may identify~$U$ with~$\iota^\sym{2}\circ U\circ(\iota^\sym{1})^{-1}\colon L^2(\mu^1)\rar L^2(\mu^2)$.

\iref{i:p:IntertwiningDirInt:1}~is~\cite[\S{II.2.3}, Example, p.~182]{Dix81}.

\iref{i:p:IntertwiningDirInt:2} Both implications are a consequence of the definition of the order structure on the direct integral of $\zeta\mapsto L^2(\mu_\zeta)$, together with Proposition~\ref{p:DIntL^2}~\iref{i:p:DIntL^2:1}. 

\iref{i:p:IntertwiningDirInt:3} Let~$t>0$ be fixed.
Assume first that~$U_\zeta$ intertwines~$T^\sym{1}_{\zeta, t}$,~$T^\sym{2}_{\zeta,t}$.
By~\cite[\S{II.2.3}, Prop.~3, p.~182]{Dix81},
\begin{align*}
U \circ T^\sym{1}_t=&\ \dint[]{Z} U_\zeta \diff\nu(\zeta) \circ \dint[]{Z} T^\sym{1}_{\zeta,t} \diff\nu(\zeta)=  \dint[]{Z} U_\zeta \circ T^\sym{1}_{\zeta,t} \diff\nu(\zeta)
\\
=&\ \dint[]{Z} T^\sym{2}_{\zeta,t} \circ U_\zeta \diff\nu(\zeta)=T^\sym{2}_t\circ U \fstop
\end{align*}

Vice versa assume that~$U$ intertwines~$T^\sym{1}_t$ and $T^\sym{2}_t$. By Proposition~\ref{p:DirInt} \iref{i:p:DirInt3} and~\cite[\S{II.2.3}, Prop.~3, p.~182]{Dix81},
\begin{align*}
\dint[]{Z} U_\zeta \circ T^\sym{1}_{\zeta,t} \diff\nu(\zeta)=U\circ T^\sym{1}_t=T^\sym{2}_t\circ U = \dint[]{Z} T^\sym{2}_{\zeta,t} \circ U_\zeta \diff\nu(\zeta)\comma
\end{align*}
hence, by a further application of~\cite[\S{II.2.3}, Prop.~3, p.~182]{Dix81}, and by~\cite[\S{II.2.3}, Prop.~2, p.~181]{Dix81},
\begin{align*}
0=&\ \norm{U\circ T^\sym{1}_t-T^\sym{2}_t\circ U}_\op =\norm{\dint[]{Z} \ttonde{U_\zeta \circ T^\sym{1}_{\zeta,t} - T^\sym{2}_{\zeta,t}\circ U_\zeta} \diff\nu(\zeta)}_\op
\\
=&\ \nu\text{-}\esssup_\zeta \norm{U_\zeta \circ T^\sym{1}_{\zeta,t} - T^\sym{2}_{\zeta,t}\circ U_\zeta}_\op\comma
\end{align*}
which concludes the proof.
\qed\end{proof}

\section{Intertwining of ergodic decompositions}

In this section we describe the structure of order isomorphisms intertwining Dirichlet forms that are not necessarily irreducible. Informally, one might expect that one can simply decompose the Dirichlet forms into their ``connected components'' and apply the known result for irreducible Dirichlet forms. 
However, it is technically non-trivial to make this notion of ``connected components'' rigorous for abstract Dirichlet forms.

We rely here on the notion of ergodic decompositions introduced in \cite{LzDS20}. We will first show that these ergodic decompositions are essentially unique (Theorem \ref{t:Uniqueness}) and that intertwining order isomorphisms carry ergodic decompositions to ergodic decompositions (Theorem \ref{t:Main}). As a consequence, whenever the Dirichlet forms in question admit ergodic decompositions, intertwining order isomorphisms act componentwise as anticipated by the informal discussion above (Corollary \ref{c:Isomorphism}).

\subsection{Ergodic decompositions}
Let us first introduce a rigorous definition for the ergodic decomposition of a regular Dirichlet form.
\begin{definition}\label{d:ErgodicDec}
Let~$(X,\T,\mcX,\mu)$ be a locally compact Polish Radon-measure space, and~$(E,\dom{E})$ be a regular Dirichlet form on~$L^2(\mu)$.
We say that~$(E,\dom{E})$ admits an \emph{ergodic decomposition~$\zeta\mapsto (E_\zeta,\dom{E_\zeta})$ indexed by~$(Z,\mcZ,\nu)$} if there exist
\begin{enumerate}[$(i)$]
\item\label{i:d:ErgodicDec:1} a separable countably generated  probability space~$(Z,\mcZ,\nu)$ and a measurable map~$s\colon (X,\mcX)\rar (Z,\mcZ)$;
\item\label{i:d:ErgodicDec:2} a pseudo-disintegration~$\seq{\mu_\zeta}_{\zeta\in Z}$ of~$\mu$ over~$\nu$, strongly consistent with~$s$ such that if~$s^{-1}(\zeta)$ is endowed with the subspace topology and the trace $\sigma$-algebra inherited by~$(X,\T,\mcX^\mu)$, then~$(s^{-1}(\zeta),\hat\mu_\zeta)$ is a Radon-measure space for~$\nu$-a.e.~$\zeta\in Z$;
\item\label{i:d:ErgodicDec:3} a $\nu$-measurable field~$\zeta\mapsto (E_\zeta,\dom{E_\zeta})$ of regular irreducible Dirichlet forms $(E_\zeta,\dom{E_\zeta})$ on~$L^2(\mu_\zeta)$;
\end{enumerate}
such that
\begin{align*}
L^2(\mu)=\dint[]{Z} L^2(\mu_\zeta)\diff\nu(\zeta) \qquad \text{and} \qquad E=\dint[]{Z} E_\zeta \diff\nu(\zeta)
\end{align*}
\end{definition}
\purple{as a direct integral of Dirichlet forms in the sense of Definition~\ref{d:DirIntE}.}

\begin{remark} Since the pseudo-disintegration~$\seq{\mu_\zeta}_{\zeta\in Z}$ in Definition~\ref{d:ErgodicDec} \iref{i:d:ErgodicDec:2} is $s$-separated, the map~$s$ is implicitly always assumed to be surjective.
\end{remark}

\begin{remark}\label{r:LiLin}
\purple{
The existence of ergodic decompositions has been shown:
\begin{enumerate}[$(a)$]
\item in~\cite{LzDS20} under the additional assumptions that either~$\mu$ be a finite measure, or that the form~$(E,\dom{E})$ be strongly local and admitting a carr\'e du champ operator;
\item\label{i:Kuwae} by K.~Kuwae in~\cite{Kuw21} under the additional assumption that the transition probabilities of the Markov process properly associated to the Dirichlet form be absolutely continuous w.r.t.\ the reference measure~$\mu$.
\end{enumerate}
It is one further result of~\cite{Kuw21} that, under the same assumption as in~\ref{i:Kuwae} above, the set~$Z$ indexing the ergodic decomposition is in fact \emph{at most countable}, which greatly simplifies the discussion.
This type of ergodic decomposition however rules out some interesting examples, in particular from infinite-dimensional analysis, see~\cite[\S3.4]{LzDS20} and references therein.
}

\purple{Let us also point out that, again in the case of~\ref{i:Kuwae} above, our main result (i.e., the extension of the results in~\cite{LenSchWir18} to the non-irreducible case) has already been obtained by L.~Li and H.~Lin in~\cite[Thm.~4.6]{LiLin22}.}
\end{remark}

Ergodic decompositions of Dirichlet forms are \emph{essentially projectively unique}, in the sense of the next theorem, expanding the scope of all the projective uniqueness results in~\cite{LzDS20}.

\begin{theorem}\label{t:Uniqueness}
Let~$(X,\mcX,\mu)$ be a locally compact Polish Radon-measure space, and~$(E,\dom{E})$ be a regular Dirichlet form on~$L^2(\mu)$. Further assume that $(E,\dom{E})$ admits
\begin{itemize}
\item an ergodic decomposition~$\zeta\mapsto (E^1_\zeta,\dom{E^1_\zeta})$ indexed by~$(Z^1,\mcZ^1,\nu^1)$;
\item an ergodic decomposition~$\eta\mapsto (E^2_\eta,\dom{E^2_\eta})$ indexed by~$(Z^2,\mcZ^2,\nu^2)$.
\end{itemize}

Then, there exists an almost isomorphism~$\varrho\colon (Z^1,\mcZ^1,\mcN^{\nu^1})\rar (Z^2,\mcZ^2,\mcN^{\nu^2})$ and a~$\mcZ^2$-measurable function~$h\colon Z^2\rar (0,\infty)$ so that~$\varrho_\pfwd \nu^1\sim \nu^2$ and~$E^1_{\zeta} = h\ttonde{\varrho(\zeta)} E^2_{\varrho(\zeta)}$ for $\nu^1$-a.e.~$\zeta\in Z^1$.
\end{theorem}

\begin{proof}
Let~$\mcX_0$ be the family of $\mu$-measurable $E$-invariant subsets of~$X$, and note that~$\mcX_0$ is a $\sigma$-subalgebra of~$\mcX^\mu$, e.g.~\cite[Lem.~1.6.1, p.~53]{FukOshTak11}.
Let~$\mu_0$ be the restriction of~$\hat\mu$ to~$(X,\mcX_0)$.
By our Assumption~\ref{ass:Main},~$\mcX$ is countably generated, thus~$\mcX_0$ is $\mu_0$-essentially countably generated by $\mcX^*\eqdef \mcX\cap \mcX_0$.

Further let~$\mcC$ be a special standard core for~$(E,\dom{E})$ witnessing the compatibility of the direct integral representation, i.e.\ so that~$\mcC=\mcA$ as in Definition~\ref{d:Compat}.
Note that for every~$u\in\mcC$ we may choose~$u_\zeta\equiv u$ as a $\mu_\zeta$-representative of the diagonal embedding~$\delta(u)_\zeta$.

\paragraph{Step~1: Invariant sets}
For simplicity of notation, in this step let~$Z$ denote either~$Z^1$ or~$Z^2$, and analogously for all the other symbols.

\paragraph{Claim I: For every~$B\in\mcZ$ the set~$A\eqdef s^{-1}(B)$ is $E$-invariant}
Since~$\seq{\mu_\zeta}_{\zeta\in Z}$ is $s$-separated,~$\seq{X_\zeta}_{\zeta\in Z}$, with~$X_\zeta\eqdef s^{-1}(\zeta)$, is a separating family in the sense of Definition~\ref{d:Disint}.
The multiplication operator~$M_{X_\zeta}\eqdef M_{\car_{X_\zeta}}\colon L^2(\mu_\zeta)\to L^2(\mu_\zeta)$ satisfies~$M_{X_\zeta}=\id_{L^2(\mu_\zeta)}$, hence~$T_{\zeta,t} M_{X_\zeta}=M_{X_\zeta} T_{\zeta,t}$.
As a consequence,
\begin{align*}
M_A \circ T_t=&\ \dint[]{B} M_{X_\zeta} \diff\nu(\zeta) \circ \dint[]{Z} T_{\zeta, t} \diff\nu(\zeta) = \dint[]{Z} M_{X_\zeta} \circ T_{\zeta,t}\diff\nu(\zeta)
\\
=&\ \dint[]{Z} T_{\zeta,t}\circ M_{X_\zeta}\diff\nu(\zeta)= \dint[]{Z} T_{\zeta, t} \diff\nu(\zeta) \circ \dint[]{B} M_{X_\zeta} \diff\nu(\zeta)
\\
=&\ T_t\circ M_A \comma
\end{align*}
that is,~$A\eqdef s^{-1}(B)$ is $E$-invariant.

\paragraph{Claim II: for every $E$-invariant~$A\in\mcX_0$ there exists~$B\in\mcZ$ with $\mu_0(A\triangle s^{-1}(B))=0$} 
It suffices to show the statement for~$A\in\mcX^*$.
Since~$(Z,\mcZ)$ is separable countably generated, it is countably separated, hence by~\cite[Prop.~12.1(iii)]{Kec95} there exists a separable metrizable topology~$\T_Z$ on~$Z$ so that~$\mcZ$ is the Borel $\sigma$-algebra generated by~$\T_Z$.
As a consequence, $s\colon X\to Z$ is Borel, and therefore~$s(A)$ is analytic by~\cite[423G(b)]{Fre00}, thus $\mcZ$-universally measurable by~\cite[434D(c)]{Fre00}, hence finally $\nu$-measurable.
In particular, there exist~$B_0\subset s(A)\subset B_1$ with~$B_0,B_1\in \mcZ$ and~$\nu(B_1\setminus B_0)=0$.
By the previous claim, we have that~$A_0\eqdef s^{-1}(B_0)$, $A_1\eqdef s^{-1}(B_1)\in \mcX_0$.
The conclusion follows since~$A_0\subset A\subset A_1$ by definition, and since~$\mu_0(A_1\setminus A_0)=0$ combining~\eqref{eq:Disint:0} with~$\nu(B_1\setminus B_0)=0$.

\paragraph{Claim III: $s^*\mcZ\eqdef\set{s^{-1}(B): B\in\mcZ}$ is $\mu_0$-essentially countably generating~$\mcX_0$} 
It suffices to combine the previous two claims with the fact that~$\mcZ$ is countably generated by assumption.

\paragraph{Claim IV: The measure space~$(Z,\mcZ,\nu)$ is perfect}
We shall use Lemma~\ref{l:Perfect} in numerous instances, without explicit mention.
As a Radon space~$(X,\T,\mcX^\mu,\mu)$ is perfect. Since it is $\sigma$-finite, there exists a probability measure~$\lambda'$ equivalent to~$\hat\mu$ on~$\mcX^\mu$, and the space~$(X,\mcX^{\lambda'},\lambda')$ is perfect as well.
Let~$\lambda\eqdef \lambda'\mrestr{s^*\mcZ}$ and~$\sigma\eqdef s_\pfwd\lambda$, and note that the restricted measure space~$(X, s^*\mcZ,\lambda)$ is perfect. Thus $(Z,\mcZ,\sigma)$ is perfect, too.

Now, let~$B\in\mcZ$ with~$\nu B=0$. Then,~$s^{-1}(B)$ is $\mu$-negligible by~\eqref{eq:Disint:0}.
Since~$\lambda\sim\mu$ on~$s^*\mcZ$, the set~$s^{-1}(B)$ is as well $\lambda$-negligible, hence~$B$ is as well $\sigma$-negligible, which shows that~$\sigma\ll\nu$.
Vice versa let~$B\in\mcZ$ with~$\sigma B=0$. Then,~$s^{-1}(B)$ is $\lambda$-negligible, and therefore $\mu$-negligible.
Since~$\mu_\zeta$ is non-zero for every~$\zeta\in Z$ by definition of pseudo-disintegration, we conclude again from~\eqref{eq:Disint:0} that~$B$ is as well $\nu$-negligible, which shows $\nu\ll\sigma$.
Thus~$\nu\sim\sigma$.

Finally, since both~$\nu$ and~$\sigma$ are finite measures,~$\nu$ has a density w.r.t.~$\sigma$ by the Radon--Nikodym Theorem, and therefore~$(Z,\mcZ,\nu)$ is perfect since so is~$(Z,\mcZ,\sigma)$.

\paragraph{Step~2: Almost isomorphism}
Throughout this step, let~$i=1,2$.
By construction, the measure algebra~$(\mfZ_i, \bar\sigma_i)$ of~$(Z,\mcZ^i,\sigma_i)$ is isomorphic to the measure algebra of~$(X, s_i^*\mcZ^i,\lambda^i)$.
Furthermore, by \emph{Step~1: Claim III} the $\sigma$-algebra $\mcZ^i$ is $\mu_0$-essentially countably generating the algebra~$\mcX_0$ of $E$-invariant sets.
Thus,~$\mcZ^i$ is as well $\lambda_0\eqdef \lambda'\mrestr{\mcX^*}$-essentially countably generating~$\mcX_0$, and therefore
\begin{align*}
\ttonde{X,(s_1^*\mcZ^1)^{\lambda_1},\hat\lambda_1}=\ttonde{X,(\mcX^*)^{\lambda_0},\hat\lambda_0}=\ttonde{X,(s_2^*\mcZ^2)^{\lambda_2},\hat\lambda_2}\fstop
\end{align*}
As a consequence, since the measure algebra of a measure space coincides with that of its completion, the measure algebras~$(\mfZ^1, \bar\sigma^1)$ and~$(\mfZ^2, \bar\sigma^2)$ are isomorphic.
By all of the above, the spaces~$(Z^i,(\mcZ^i)^{\sigma_i},\hat\sigma_i)$ are countably separated perfect complete probability spaces, with isomorphic measure algebras.
Therefore, the corresponding measure spaces~$(Z^i,\mcZ^i,\sigma^i)$ are almost isomorphic by~\cite[344C]{Fre00}, via some almost isomorphism~$\varrho\colon Z^1\to Z^2$ satisfying~$\varrho_\pfwd \sigma^1=\sigma^2$.
Since~$\sigma^i\sim \nu^i$ by the proof of \emph{Step 1 Claim IV}, we conclude that~$\varrho_\pfwd \nu^1\sim \nu^2$.

\paragraph{Step 3: Forms identifications}
Since~$\varrho_\pfwd \nu^1\sim \nu^2$, there exists a measurable $h\colon (Z^2,\mcZ^2)\to(0,\infty)$ satisfying~$h \cdot \varrho_\pfwd\nu^1=\nu^2$.
Now
\begin{equation}\label{eq:t:Uniqueness:1}
\begin{aligned}
\mu\ttonde{A\cap (s^2)^{-1}(B^2)}=&\ \int_{B^2} \mu^2_\eta A \diff\nu^2(\eta)
\\
=&\ \int_{\varrho^{-1}(B^2)} \mu^2_{\varrho(\zeta)} A \cdot h(\varrho(\zeta))^{-1}\diff\nu^1(\zeta)
\end{aligned}
\comma\quad A\in\mcX\comma B^2\in\mcZ^2\comma
\end{equation}
which shows that~$\tseq{h(\varrho(\zeta))^{-1} \cdot \mu^2_{\varrho(\zeta)}}_{\zeta\in Z^1}$ is a pseudo-disintegration of~$\mu$ over~$\nu^1$, strongly consistent with~$\varrho^{-1}\circ s^2$.

Furthermore, exchanging the roles of~$\nu^1$ and $\nu^2$ in \emph{Step~2} we have that~$\varrho_\pfwd \nu^1\sim\nu^2$ and~$\varrho^{-1}_\pfwd\nu^2\sim\nu^1$.
Thus,~$f\colon Z^2\to \R$ is $\nu^2$-measurable if and only if~$f\circ\varrho\colon Z^1\to\R$ is $\nu^1$-measurable.
As a consequence, the function~$h\circ \varrho\colon Z^1\to (0,\infty)$ is $\nu^1$-measurable and, since~$\eta\mapsto (E^2_\eta,\dom{E^2_\eta})$ is a $\nu^2$-measurable field of quadratic forms, for every~$u,v\in\mcC$ the map
\begin{align*}
\zeta\longmapsto E^2_{\varrho(\zeta),1}(u,v) = E^2_{\varrho(\zeta)}(u_{\varrho(\zeta)},v_{\varrho(\zeta)})+\scalar{u_{\varrho(\zeta)}}{v_{\varrho(\zeta)}}_{L^2(\mu^2_{\varrho(\zeta)})}
\end{align*}
is $\nu^1$-measurable, which shows that $\zeta\mapsto \ttonde{h(\rho(\zeta)) E^2_{\varrho(\zeta)},\dom{E^2_{\varrho(\zeta)}}}$ is a $\nu^1$-measurable field of quadratic forms on~$L^2(\mu^2_{\varrho(\zeta)})$.
The corresponding direct integral form~$\ttonde{\widetilde{E},\dom{\widetilde{E}}}$ is compatible with the disintegration in~\eqref{eq:t:Uniqueness:1} in the sense of Definition~\ref{d:Compat} with underlying space~$S_\mcC$, since~$\mcC$ is a core for
\begin{equation*}
\ttonde{h(\rho(\zeta)) E^2_{\varrho(\zeta)},\dom{E^2_{\varrho(\zeta)}}}
\end{equation*}
for $\nu^1$-a.e.~$\zeta\in Z^1$.
Finally, let us show that the form~$\ttonde{\widetilde{E},\dom{\widetilde{E}}}$ coincides with~$(E,\dom{E})$.
It suffices to note that, for all~$u,v\in\mcC$ we have
\begin{align*}
\widetilde{E}(u,v) =&\ \int_{Z^1} h(\varrho(\zeta)) E^2_{\varrho(\zeta)}(u_{\varrho(\zeta)}, v_{\varrho(\zeta)})\diff\nu^1(\zeta)
\\
=&\ \int_{Z^2} h(\eta) E^2_\eta(u_\eta,v_\eta)\diff\varrho_\pfwd \nu^1= \int_{Z^2} E^2_\eta(u_\eta, v_\eta)\diff\nu^2(\eta)=E(u,v) \fstop \qquad \qed
\end{align*}
\end{proof}

\subsection{Intertwining}
We are now ready to state our main result.

\begin{theorem}\label{t:Main}
For~$i=1,2$, let~$(X^i,\T^i,\mcX^i,\mu^i)$ be a locally compact Polish Radon-measure space, and~$(E^i,\dom{E^i})$ be a regular Dirichlet form on~$L^2(\mu^i)$ with semigroup~$\seq{T^\sym{i}_t}_{t\geq 0}$.
Suppose further that
\begin{enumerate}[$(a)$]
\item $(E^1,\dom{E^1})$ admits ergodic decomposition~$\zeta\mapsto (E^1_\zeta,\dom{E^1_\zeta})$ indexed by $(Z,\mcZ,\nu)$;
\item there exists a \emph{unitary} order isomorphism~$U\colon L^2(\mu^1)\rar L^2(\mu^2)$ intertwining~$\seq{T^\sym{1}_t}_{t\geq 0}$,~$\seq{T^\sym{2}_t}_{t\geq 0}$, represented by a scaling~$h$ and a transformation~$\tau$ in the sense of Proposition~\ref{p:Lamperti}.
\end{enumerate}

Finally, set
\begin{align}\label{eq:t:Main:0}
\mu^2_\zeta=\frac{\diff\mu^2}{\diff (\tau^{-1})_\pfwd\mu^1}(\tau^{-1})_\pfwd \mu^1_\zeta \qquad \forallae{\nu}\zeta\in Z\fstop
\end{align}

Then,
\begin{enumerate}[$(i)$]
\item\label{i:t:Main:0} $\tseq{\mu^2_\zeta}_{\zeta\in Z}$ is a separated pseudo-disintegration of~$\mu^2$ over~$(Z,\mcZ,\nu)$;

\item\label{i:t:Main:1} $U\colon L^2(\mu^1)\rar L^2(\mu^2)$ is decomposable over~$(Z,\mcZ,\nu)$ and represented by a $\nu$-measurable field of order isomorphisms~$U_\zeta\colon L^2(\mu^1_\zeta)\rar L^2(\mu^2_\zeta)$.
\end{enumerate}

For~$\zeta\in Z$ and~$t>0$ further define an operator~$T^\sym{2}_{\zeta,t}\colon L^2(\mu^2_\zeta)\rar L^2(\mu^2_\zeta)$ by
\begin{align}\label{eq:T2zetaDef}
T^\sym{2}_{\zeta, t} \circ U_\zeta = U_\zeta \circ T^\sym{1}_{\zeta, t} \fstop
\end{align}
Then,

\begin{enumerate}[$(i)$]\setcounter{enumi}{2}
\item\label{i:t:Main:1.5} $\tseq{T^\sym{2}_{\zeta,t}}_{t\geq 0}$ is a sub-Markovian strongly continuous contraction semigroup on~$L^2(\mu^2_\zeta)$ for $\nu$-a.e.~$\zeta\in Z$;

\item\label{i:t:Main:2} $(E^2,\dom{E^2})$ admits an ergodic decomposition~$\zeta\mapsto (E^2_\zeta, \dom{E^2_\zeta})$ indexed by~$(Z,\mcZ,\nu)$, where the Dirichlet form~$(E^2_\zeta,\dom{E^2_\zeta})$ is the image form of $(E^1_\zeta,\dom{E^1_\zeta})$ via~$U_\zeta$ for $\nu$-a.e.~$\zeta\in Z$ and the associated semigroup is~$\tseq{T^\sym{2}_{\zeta,t}}_{t\geq 0}$.
\end{enumerate}
\end{theorem}

\begin{proof}
Up to (re-)defining~$\tau$, resp.~$\tau^{-1}$, on a $\mu^1$-, resp.~$\mu^2$-, negligible set, we can assume without loss of generality that~$\tau$ is defined everywhere. Note however that~$\tau$ is neither surjective nor injective.

\paragraph{Step 1: Disintegrations} By Lemma~\ref{l:Adjoint}, the definition of~$\tseq{\mu^2_\zeta}_{\zeta\in Z}$ in~\eqref{eq:t:Main:0} is well-posed.
By definition of pseudo-disintegration,~$\zeta\mapsto \mu^1_\zeta f$ is $\nu$-measurable for every $\mcX^1$-measurable real-valued~$f$ on~$X^1$. As a consequence, the map
\begin{align*}
\zeta \mapsto \mu^2_\zeta A = \int_A \frac{\diff\mu^2}{\diff (\tau^{-1})_\pfwd \mu^1} \diff (\tau^{-1})_\pfwd \mu^1_\zeta \comma \qquad A\in \mcX^2\comma
\end{align*}
is as well $\nu$-measurable, by $\mcX^2/\mcX^1$-measurability of~$\tau$.
Furthermore, by~\eqref{eq:Disint:0} for~$\tseq{\mu^1_\zeta}_{\zeta\in Z}$,
\begin{align*}
\int_Z \int_A \frac{\diff\mu^2}{\diff (\tau^{-1})_\pfwd \mu^1} \diff (\tau^{-1})_\pfwd \mu^1_\zeta \diff\nu(\zeta)=& \int_Z\int_A \frac{\diff\mu^2}{\diff (\tau^{-1})_\pfwd \mu^1}\circ\tau^{-1} \diff\mu^1_\zeta \diff\nu(\zeta)
\\
=&\ \int_A \frac{\diff\mu^2}{\diff (\tau^{-1})_\pfwd \mu^1}\circ\tau^{-1} \diff\mu^1\\ 
=& \int_A \frac{\diff\mu^2}{\diff (\tau^{-1})_\pfwd \mu^1} \diff(\tau^{-1})_\pfwd\mu^1
\\
=&\ \mu^2 A \comma
\end{align*}
that is, the family~$\tseq{\mu^2_\zeta}_{\zeta}$ is a pseudo-disintegration of~$\mu^2$ over~$(Z,\mcZ,\nu)$.
Since~$\tau^{-1}$ is $\mu^1$-a.e.\ injective, it is as well~$\mu^1_\zeta$-a.e.\ injective for $\nu$-a.e~$\zeta\in Z$.
As a consequence, the pseudo-disintegration~$\tseq{\mu^2_\zeta}_{\zeta}$ is as well separated, which concludes a proof of~\iref{i:t:Main:0}.

\paragraph{Step~2: Direct-integral representations of $L^2$-spaces} Since~$\tseq{\mu^1_\zeta}_{\zeta\in Z}$ is separated by assumption and since~$\tseq{\mu^2_\zeta}_{\zeta\in Z}$ is separated by~\iref{i:t:Main:0}, Proposition~\ref{p:DIntL^2} yields the direct-integral representations
\begin{align*}
L^2(\mu^i)=\dint[S^i]{Z} L^2(\mu^i_\zeta)\diff\nu(\zeta) \comma \qquad i=1,2\comma
\end{align*}
for every space of $\nu$-measurable vector fields~$S^i=S_{\mcA^i}$ induced by~$\mcA^i$ as in~\eqref{eq:AssA}.
In the following, we shall always choose~$\mcA^1\eqdef \mcC$ to be a special standard core for the (regular) form~$(E^1,\dom{E^1})$, and set~\purple{$\mcA^2\eqdef \set{h\cdot f\circ\tau : f \in\mcA^1}$.
In light of Proposition~\ref{p:IntertwiningDirichlet} we have~$\mcA^2\subset \dom{E^2}$.
}
Having fixed the spaces~$S^i$ throughout the proof, we omit them from the notation.

As a consequence, \emph{if}~$U\colon L^2(\mu^1)\rar L^2(\mu^2)$ is decomposable, and represented by a $\nu$-measurable field of operators~$\zeta\mapsto U_\zeta$, then~$U_\zeta\colon L^2(\mu^1_\zeta)\rar L^2(\mu^2_\zeta)$ for $\nu$-a.e.~$\zeta\in Z$.

\paragraph{Step~3: Decomposability of~$U$}
For every~$A\in(\mcX^1)^{\mu^1}$, and every~$f\in L^2(\mu^1)$,
\begin{align}\label{eq:t:Main:1}
\begin{aligned}
U M_A f =&\ h\cdot \tau^*(\car_A f) = h \cdot \tau^*\car_A \cdot \tau^*f 
\\
=&\ \tau^*\car_A \cdot h \cdot \tau^* f= \tau^*\car_A \cdot Uf=M_{\tau^{-1}(A)} U f \fstop
\end{aligned}
\end{align}

Now, let~$A$ be $E^1$-invariant.
Combining~\eqref{eq:t:Main:1} and the definition of invariant set yields
\begin{align*}
M_{\tau^{-1}(A)} \circ T^\sym{2}_t \circ U=&\ M_{\tau^{-1}(A)} \circ U \circ T^\sym{1}_t\\
=&\ U \circ M_A \circ T^\sym{1}_t\\
=&\ U \circ T^\sym{1}_t \circ M_A\\
=&\ T^\sym{2}_t \circ U \circ M_A
\\
=&\ T^\sym{2}_t \circ M_{\tau^{-1}(A)} \circ U \comma
\end{align*}
whence, pre-composing with~$U^{-1}$ on both sides yields that~$\tau^{-1}(A)$ is $E^2$-invariant.

Arguing as in \emph{Step 3} in the proof of~\cite[Thm.~3.4]{LzDS20}, it follows that~$U$ commutes with all the operators in $\mcB(L^2(\mu^1))$ diagonalizable over~$Z$ in the sense of~\cite[\S{II.2.4}, Dfn.~3, p.~185]{Dix81}.
By the characterization of decomposable operators via diagonalizable operators~\cite[\S{II.2.5}, Thm.~1, p.~187]{Dix81}, the operator~$U$ is decomposable, and represented by a $\nu$-measurable field of bounded operators~$\zeta\mapsto U_\zeta\colon L^2(\mu^1_\zeta)\rar L^2(\mu^2_\zeta)$.
For $\nu$-a.e.~$\zeta\in Z$, the operator~$U_\zeta$ is additionally an order isomorphism, by the forward implication in Proposition~\ref{p:IntertwiningDirInt}. This concludes the proof of~\iref{i:t:Main:1}.

\paragraph{Step 4: Image forms} 
Since~$U_\zeta$ is unitary and $\tseq{T^\sym{1}_{\zeta,t}}_{t\geq 0}$ is a strongly continuous contraction semigroup on~$L^2(\mu^1_\zeta)$ for $\nu$-a.e.~$\zeta\in Z$, the same holds for~$\tseq{T^\sym{2}_{\zeta,t}}_{t\geq 0}$.
For~$t>0$, let
\begin{align*}
\tilde T^\sym{2}_t\eqdef \dint[]{Z} T^\sym{2}_{\zeta,t}\diff\nu(\zeta)\fstop
\end{align*}
By \emph{Step~2} above,~$\tilde T^\sym{2}_t$ is a bounded operator in~$\mcB(L^2(\mu^2))$. By the decomposability of~$U$ shown in \emph{Step~4} and by definition~\eqref{eq:T2zetaDef} of~$T^\sym{2}_{\zeta,t}$, we may apply the reverse implication in Proposition~\ref{p:IntertwiningDirInt} \iref{i:p:IntertwiningDirInt:3} to obtain
\begin{align*}
\tilde T^\sym{2}_t\circ U= U\circ T^\sym{1}_t\comma \qquad t>0\comma
\end{align*}
and conclude that~$\tilde T^\sym{2}_t=T^\sym{2}_t$ for each~$t>0$, since~$U$ is unitary.

\purple{For~$\nu$-a.e.~$\zeta\in Z$, let~$\ttonde{E^2_\zeta,\dom{E^2_\zeta}}$ be the quadratic form on~$L^2(\mu^2_\zeta)$ corresponding to~$\tseq{T^\sym{2}_{\zeta,t}}_{t\geq 0}$.}

\paragraph{Step~5: Direct-integral representation}
The direct integral representation of $(E^2,\dom{E^2})$ \purple{as a direct integral of Dirichlet forms} follows from that of~$\seq{T^\sym{2}_t}_{t\geq 0}$ by Proposition~\ref{p:DirInt}, provided we show
\begin{enumerate*}[$(a)$]
\item\label{i:Step5.1} the compatibility of~$(E^2,\dom{E^2})$ with the representation by~$\zeta\mapsto E^2_\zeta$ in the sense of Definition~\ref{d:Compat}; \purple{and that}
\item\label{i:Step5.2}\purple{$\ttonde{E^2_\zeta,\dom{E^2_\zeta}}$ is a Dirichlet form for $\nu$-a.e.~$\zeta\in Z$.}
\end{enumerate*}
In light of \emph{Step~2}, \purple{in order to show~\ref{i:Step5.1}} it suffices to show that~$\mcA^2$ is contained in~$\dom{E^2}$, which follows from Proposition~\ref{p:IntertwiningDirichlet}.
\purple{In order to show~\ref{i:Step5.2}, note that, since~$\seq{T^\sym{2}_t}_{t\geq 0}$ is sub-Markovian, it follows by Prop.~\ref{p:DirIntE} that~$\tseq{T^\sym{2}_{\zeta,t}}_{t\geq 0}$ is sub-Markovian for $\nu$-a.e.~$\zeta\in Z$.
}
\qed\end{proof}

\begin{corollary}\label{c:Isomorphism}
For~$i=1,2$ let~$(X^i,\T^i,\mcX^i,\mu^i)$ be a locally compact Polish Radon-measure space, and~$(E^i,\dom{E^i})$ be a regular Dirichlet form on~$L^2(\mu^i)$ with semigroup~$\seq{T^\sym{i}_t}_{t\geq 0}$, each admitting ergodic decomposition~$\zeta\mapsto \ttonde{E^i_\zeta,\dom{E^i_\zeta}}$ indexed by~$(Z^i,\mcZ^i,\nu^i)$.
Suppose further that there exists a unitary order isomorphism~$U\colon L^2(\mu^1)\rar L^2(\mu^2)$ intertwining~$\seq{T^\sym{1}_t}_{t\geq 0}$,~$\seq{T^\sym{2}_t}_{t\geq 0}$, represented by a scaling~$h$ and a transformation~$\tau$ in the sense of Proposition~\ref{p:Lamperti}.
Then, there exists an almost isomorphism~$\varrho \colon (Z^1,\mcZ^1,\nu^1)\rar (Z^2,\mcZ^2,\nu^2)$ such that
\begin{enumerate}[$(i)$]
\item $U\colon L^2(\mu^1)\rar L^2(\mu^2)$ is decomposable over~$(Z^1,\mcZ^1,\nu^1)$ and represented by a $\nu^1$-measurable field of order isomorphisms~$U_\zeta\colon L^2(\mu^1_\zeta)\rar L^2(\mu^2_{\varrho(\zeta)})$;
\item $U^{-1}\colon L^2(\mu^2)\rar L^2(\mu^1)$ is decomposable over~$(Z^2,\mcZ^2,\nu^2)$. Furthermore, the field of order isomorphisms $(U_\zeta)^{-1}\colon L^2(\mu^2_{\varrho(\zeta)})\rar L^2(\mu^2_\zeta)$ is $\nu^2$-measur\-able, and it represents~$U^{-1}$;
\item for $\nu^1$-a.e.~$\zeta\in Z^1$, the order isomorphism~$U_\zeta$ intertwines the semigroup $\tseq{T^\sym{1}_{\zeta,t}}_{t\geq 0}$ associated to~$\ttonde{E^1_\zeta,\dom{E^1_\zeta}}$ with the semigroup~$\tseq{T^\sym{2}_{\varrho(\zeta),t}}_{t\geq 0}$ associated to $\ttonde{E^2_{\varrho(\zeta)},\dom{E^2_{\varrho(\zeta)}}}$.
\end{enumerate}
\end{corollary}

\purple{
\begin{example}
Theorem~\ref{t:Main} and Corollary~\ref{c:Isomorphism} provide an alternative proof for the construction of the ergodic decomposition of the \emph{Dirichlet--Ferguson diffusion}~\cite{LzDS17+} on the space of probability measures over a closed Riemannian manifold, obtained by ad hoc techniques in~\cite[Thm.~6.15$(ii)$]{LzDS17+}.
The ergodic decomposition of the associated form consists of uncountably many components, and therefore the result does not follow from the aforementioned work~\cite{LiLin22}, cf.~Remark~\ref{r:LiLin} above.
\end{example}
}

Theorem~\ref{t:Main} has the following natural converse, a proof of which is a straightforward application of Proposition~\ref{p:IntertwiningDirInt}

\begin{proposition}
Let~$(X,\T,\mcX,\mu)$ be a locally compact Polish Radon-measure space, $(Z,\mcZ,\nu)$ be a separable countably generated probability space, and $s\colon (X,\mcX^\mu)\to (Z,\mcZ)$ be a measurable map.
For~$i=1,2$ further let~$\ttseq{\mu^i_\zeta}_{\zeta\in Z}$ be an $s$-separated pseudo-disintegration of~$\mu$ over~$(Z,\mcZ,\nu)$, and~$\zeta\mapsto (E^i_\zeta,\dom{E^i_\zeta})$ be a $\nu$-measurable field of Dirichlet forms compatible with the disintegration.
Further assume that there exists a $\nu$-measurable field~$\zeta\mapsto U_\zeta$ of unitary order isomorphisms $U_\zeta\colon L^2(\mu_\zeta)\to L^2(\mu_\zeta)$ intertwining the semigroups of~$\ttonde{E^i_\zeta,\dom{E^i_\zeta}}$ as in~\eqref{eq:T2zetaDef}.

\noindent Then, the direct integral forms~$E^i\eqdef\dint[]{Z} E^i_{\zeta} \diff\nu(\zeta)$ are intertwined by the unitary order isomorphism $U^i\eqdef\dint[]{Z} U_{\zeta}\diff\nu(\zeta)$.
\end{proposition}

By the transfer method for quasi-regular Dirichlet spaces, e.g.~\cite{CheMaRoe94}, it is possible to extend all the previous results to the quasi-regular case.
The Definition~\ref{d:ErgodicDec} of ergodic decomposition is readily adapted by letting~$(X,\T,\mcX,\mu)$ be satisfying Assumption~\ref{ass:Main} (as opposed to: locally compact Polish).
We only spell out the adaptation of Theorem~\ref{t:Main}. The easy adaptation of Proposition~\ref{t:Uniqueness} and Corollary~\ref{c:Isomorphism} is left to the reader.

\begin{corollary}
For~$i=1,2$, let~$(X^i,\T^i,\mcX^i,\mu^i)$ be satisfying Assumption~\ref{ass:Main}, and~$(E^i,\dom{E}^i)$ be a quasi-regular Dirichlet form on~$L^2(\mu^i)$ with semigroup $\seq{T^\sym{i}_t}_{t\geq 0}$.
Let all other assumptions of Theorem~\ref{t:Main} hold. Then, all the conclusions of that Theorem hold as well.
\end{corollary}

\purple{
\begin{remark}[On quasi-homeomorphisms]
It is shown in~\cite[Thm.~3.11]{LenSchWir18} that the transformation~$\tau$ associated to any order isomorphism~$U$ intertwining two \emph{irreducible} regular Dirichlet forms has a version~$\tilde\tau$ which is additionally a quasi-homeomorphism.
This result has been extended by~L.~Li and H.~Lin in~\cite{LiLin22} to the case of up-to-countable ergodic decompositions as in Remark~\ref{r:LiLin}\ref{i:Kuwae}.
Whereas there is a reasonable expectation for this result to hold ---at least under the additional assumption that~$U$ be unitary, as in Theorem~\ref{t:Main}--- a proof seems currently beyond reach.
\end{remark}
}

\purple{
\section*{Acknowledgements}
The authors are grateful to an anonymous Reviewer for their careful reading of the manuscript and useful suggestions.
}

%KINDLY STOP ASKING FOR THIS NONSENSE
{\footnotesize
\bigskip

\noindent\textbf{Conflict of interest.} The authors declare that they have no conflict of interest.

\noindent\textbf{Data availability.} The manuscript has no associated data.
}

%\bibliographystyle{spmpsci}
%\bibliography{Bibliography.bib}  

\begin{thebibliography}{10}
\providecommand{\url}[1]{{#1}}
\providecommand{\urlprefix}{URL }
\expandafter\ifx\csname urlstyle\endcsname\relax
  \providecommand{\doi}[1]{DOI~\discretionary{}{}{}#1}\else
  \providecommand{\doi}{DOI~\discretionary{}{}{}\begingroup
  \urlstyle{rm}\Url}\fi

\bibitem{AreBietEl12}
Arendt, W., Biegert, M., ter Elst, A.F.M.: {Diffusion determines the manifold}.
\newblock {J.\ reine angew.\ Math.} \textbf{667}, 1--25 (2012).
\newblock \doi{10.1515/crelle.2011.131}

\bibitem{CheMaRoe94}
{Chen, Z.-Q.}, {Ma, Z.-M.}, {R\"ockner, M.}: {Quasi-homeomorphisms of Dirichlet
  forms}.
\newblock {Nagoya Math. J.} \textbf{136}, 1--15 (1994)

\bibitem{LzDS20}
Dello~Schiavo, L.: {Ergodic Decomposition of Dirichlet Forms via Direct
  Integrals and Applications}.
\newblock {Potential Anal.}  (2021).
\newblock \doi{10.1007/s11118-021-09951-y}.
\newblock {43 pp.}

\bibitem{LzDS17+}
Dello~Schiavo, L.: {The Dirichlet--Ferguson Diffusion on the Space of
  Probability Measures over a Closed Riemannian Manifold}.
\newblock {Ann.\ Probab.} \textbf{50}(2), 591--648 (2022).
\newblock \doi{10.1214/21-AOP1541}.
\newblock {57 pp.}

\bibitem{Dix81}
{Dixmier, J.}: {Von Neumann Algebras}.
\newblock {North-Holland} (1981)

\bibitem{Fre00}
{Fremlin, D.~H.}: {Measure Theory -- Volume I - IV, V Part I \& II}.
\newblock {Torres Fremlin (ed.)} (2000-2008)

\bibitem{FukOshTak11}
{Fukushima, M.}, {Oshima, Y.}, {Takeda, M.}: {Dirichlet forms and symmetric
  Markov processes}, \emph{{De Gruyter Studies in Mathematics}}, vol.~19,
  extended edn.
\newblock {de Gruyter} (2011)

\bibitem{HayMirYve91}
{Haydon, Richard}, {Levy, Mireille}, {Raynaud, Yves}: {Randomly normed spaces}.
\newblock {Travaux en cours}. Hermann (1991)

\bibitem{Kac66}
Kac, M.: {Can One Hear the Shape of a Drum?}
\newblock {Amer.\ Math.\ Month.} \textbf{73}(4), 1--23 (1966)

\bibitem{Kec95}
{Kechris, A.~S.}: {Classical Descriptive Set Theory}, \emph{{Graduate Texts in
  Mathematics}}, vol. 156.
\newblock {Springer-Verlag}, {New York} ({1995})

\bibitem{KelLenSchWir15}
{Keller, M.}, {Lenz, D.}, {Schmidt, M.}, {Wirth, M.}: {Diffusion determines the
  recurrent graph}.
\newblock {Adv.\ Math.} \textbf{269}, 364--398 (2015)

\bibitem{Kuw21}
Kuwae, K.: {Irreducible decomposition for Markov processes}.
\newblock {Stoch.\ Proc.\ Appl.} \textbf{140}, 339--356 (2021).
\newblock \doi{10.1016/j.spa.2021.06.012}

\bibitem{LenSchWir18}
{Lenz, D.}, {Schmidt, M.}, {Wirth, M.}: {Geometric Properties of Dirichlet
  Forms under Order Isomorphism}.
\newblock {arXiv:1801.08326}  (2018)

\bibitem{LiLin22}
Li, L., Lin, H.: {On Order Isomorphisms Intertwining Semigroups for Dirichlet
  Forms}.
\newblock {arXiv:2204.02975}  (2022)

\bibitem{MaRoe92}
{Ma, Z.-M.}, {R\"ockner, M.}: Introduction to the Theory of (Non-Symmetric)
  Dirichlet Forms.
\newblock {Graduate Studies in Mathematics}. Springer (1992)

\bibitem{Wei84}
{Weis, L.}: {On the Representation of Order Continuous Operators by Random
  Measures}.
\newblock {Trans.\ Amer.\ Math.\ Soc.} \textbf{285}(2), 535--563 (1984)

\end{thebibliography}

\end{document}